\input amstex.tex
\input amsppt.sty   
\magnification 1200
\vsize = 8.28 true in
\hsize=6.2 true in
\nologo
\NoRunningHeads        
\parskip=\medskipamount
        \lineskip=2pt\baselineskip=18pt\lineskiplimit=0pt
       
        \TagsOnRight
        \NoBlackBoxes

        \topmatter
        \title
        Energy Supercritical Nonlinear Schr\"o- \\
        dinger Equations: Quasi-periodic solutions
        \endtitle
\author
         W.-M.~Wang        \endauthor        
\address
{CNRS and D\'epartement de Math\'ematique, Universit\'e Cergy-Pontoise, 95302 Cergy-Pontoise Cedex, France}
\endaddress
        \email
{wei-min.wang\@math.cnrs.fr}
\endemail
\abstract
We construct time quasi-periodic solutions to the energy supercritical nonlinear Schr\"odinger equations on the torus 
in arbitrary dimensions. This introduces a {\it new} approach, which could have general applicability.
\endabstract

        \bigskip\bigskip
        \bigskip
        \toc
        \bigskip
        \bigskip 
        \widestnumber\head {Table of Contents}
        \head 1. Introduction and statement of the Theorem
        \endhead
        \head 2. The genericity conditions
        \endhead
        \head 3. The first step in the Newton scheme --$\qquad\qquad$
         extraction of parameters
        \endhead
        \head 4. The second step
        \endhead
        \head 5. Proof of the Theorem
        \endhead
         \endtoc
        \endtopmatter
        \vfill\eject
        \bigskip
\document
\head{\bf 1. Introduction and statement of the Theorem}\endhead
We consider the nonlinear Schr\"odinger equations (NLS) on the $d$-torus: $\Bbb T^d:=\Bbb R^d/(2\pi\Bbb Z)^d$:
$$
i\frac\partial{\partial t}u =-\Delta u+|u|^{2p}u+H(x,u,\bar u),\tag 1.1
$$
where $p\geq 1$ and $p\in\Bbb N$ is {\it arbirary}; $\Delta:=\sum^d_{k=1}\partial^2/\partial x^2_k$ is the usual Laplacian; $u$ is a function on $\Bbb R\times \Bbb T^d$ and for each given $t\in\Bbb R$, $u$ is identified with a periodic function on $\Bbb R^d$: 
$u(t,x)=u(t, x+2j\pi)$ for all $j\in\Bbb Z^d$, $\bar u$ is the complex conjugate of $u$;
$H(x,u,\bar u)$ is analytic in $(x, u,\bar u)$ and has the expansion:
$$H(x,u,\bar u)=\sum_{m=1}^\infty \alpha_m(x)|u|^{2p+2m} u,$$
where $\alpha_m$ as a function on $\Bbb R^d$ is $(2\pi)^d$ periodic and real and analytic in a strip of width $\Cal O (1)$
for all $m$. 

Using Fourier series, the solutions to the linear equation: 
$$i\frac\partial{\partial t}u+\Delta u=0\tag 1.2$$
are linear combinations of eigenfunction solutions of the form:
$$e^{-ij^2t}e^{ij\cdot x},\quad j\in\Bbb Z^d,$$
where $j^2:=|j|^2$ and $\cdot$ is the usual scalar product. Since the time frequencies are integers, the linear flow is periodic 
with $1$ being the basic frequency. 

After the addition of the nonlinear terms and for small $u$, it is natural to ask whether the linear solutions could bifurcate 
to solutions to the nonlinear equation (1.1), albeit with several frequencies -- the quasi-periodic solutions. In this paper,
we address this question using a space-time approach, which makes available a much larger phase-space.


\subheading {1.1 The space-time Fourier series}
  
To proceed, let $u^{(0)}$ be a solution with a finite number of frequencies to the linear equation in (1.2).
Fix the number of frequencies to be $b$ and write the solution as 
$$u^{(0)}(t, x)=\sum_{k=1}^b a_ke^{-ij_k^2t}e^{ij_k\cdot x}.$$
For the nonlinear construction, it is useful to add a dimension for each frequency in time and view $u^{(0)}$
as a function on $\Bbb T^b\times\Bbb T^d: =\Bbb T^{b+d}\supset\Bbb T^d$ -- this is called a {\it lift}. 
Henceforth $u^{(0)}$ adopts the form: 
$$
\aligned u^{(0)}(t, x)&=\sum_{k=1}^b a_ke^{-ij_k^2t}e^{ij_k\cdot x}\\
:&=\sum_{k=1}^b{\hat u}^{(0)}(-e_{k}, j_k)e^{-i(e_{k}\cdot\omega^{(0)}) t}e^{ij_k\cdot x}, \endaligned
$$
where $e_{k}= (0, 0, ...1, .., 0)\in\Bbb Z^b$ is a unit vector, with the only non-zero component in the $k$th direction, 
$\omega^{(0)}=\{j_k^2\}_{k=1}^{b}$ ($j_k\neq 0$) and $\hat u^{(0)}(-e_{k}, j_k)=a_k$.
Therefore $u^{(0)}$ has Fourier support
$$\text{supp } {\hat u}^{(0)}=\{(-e_{k}, j_k), k=1,...,b\}\subset\Bbb Z^{b+d},\tag 1.3$$
where $j_k\neq j_{k'}$ if $k\neq k'$.

For the nonlinear equation (1.1), we seek quasi-periodic solutions with $b$ frequencies in the form of a space-time Fourier series:
$$
u(t, x)=\sum_{(n,j)\in\Bbb Z^{b+d}}\frak a(n, j)e^{in\cdot\omega t}e^{ij\cdot x}. \qquad \tag 1.4
$$
We iteratively determine $\frak a$ together with the frequency $\omega\in\Bbb R^b$ . This is the well-known amplitude-frequency
modulation fundamental to nonlinear equations. We note that when the equation is linear, the frequency $\omega$ is independent of $u$
and (1.4) reduces 
to the usual Fourier series. For example for the linear solution $u^{(0)}$, 
the frequencies are {\it fixed} at $\omega=\omega^{(0)}=\{j_k^2\}_{k=1}^b\in\Bbb R^b$, which are eigenvalues of the Laplacian.

\noindent{\it Remark.} The space-time Fourier series can be viewed as resulting from the embedding: 
$$\{1, 2, ..., b\}\hookrightarrow\Bbb Z^b.$$ The ambient space $\Bbb Z^b$ is precisely the Fourier dual of $\Bbb T^b$.

In the Fourier space $\Bbb Z^{b+d}$ , the support of the solution $u$ in the form (1.4) to the linear equation (1.2) 
is the characteristics $\Cal C^+$: 
$$\Cal C^+=\{(n,j)\in\Bbb Z^{b+d}|n\cdot\omega^{(0)}+j^2=0\}.\tag 1.5$$
The support of the complex conjugate $\bar u$ is 
the characteristics $\Cal C^-$: 
$$\Cal C^-=\{(n,j)\in\Bbb Z^{b+d}|-n\cdot\omega^{(0)}+j^2=0\}.\tag 1.6$$
It is convenient to define the bi-characteristics $\Cal C$ as 
$$\Cal C:=\Cal C^+\oplus \Cal C^-\subset\Bbb Z^{b+d}\oplus \Bbb Z^{b+d}\sim\Bbb Z^{b+d}\times\Bbb Z_2.\tag 1.7$$ 
We also call $\Cal C$, the singular set, as it  can possibly contribute to the non-invertibility of a linearized operator in the Newton scheme
that we shall use starting in sect. 3.  

We consider $\Cal C$ as the restriction to $\Bbb Z^{b+d}\times \Bbb Z_2$ of the corresponding 
paraboloids on $\Bbb R^{b+d}\times \Bbb Z_2$, obtained by considering (1.5, 1.6) on $\Bbb R^{b+d}$ instead of on $\Bbb Z^{b+d}$. 
In this sense, we say that $\Cal C$ is a manifold of singularities and not just isolated
points. Moreover since $\omega^{(0)}$ is an integer vector, $\Cal C$ not only lacks convexity but also has null-directions in $n$,
defined to be the set $\{n\in\Bbb Z^b|n\cdot \omega^{(0)}=0\}$,  which is an infinite set.

Assume that the linear solution $u^{(0)}$:
$$u^{(0)}(t, x)=\sum_{k=1}^b a_ke^{-ij_k^2t}e^{ij_k\cdot x},\tag 1.8$$
is generic, satisfying the genericity conditions (Gi-iv) in sect. 2.  
They pertain entirely to the spatial Fourier support of 
$u^{(0)}$: $\{j_k\}_{k=1}^b\in (\Bbb Z^d)^b\subset (\Bbb R^d)^b $ and are determined by the $|u|^{2p}u$ term in (1.1) only. 

Assume that the time frequency satisfies: 
$$n\cdot\omega^{(0)}\neq 0,\tag $\flat$ $$
for $n\in [-N, N]^b\backslash\{(0)\}$,
where $N$ is assumed to be large.
Assume that the dimension of the torus $b$ is large satisfying $$b>C_pd.\tag$\flat\flat$ $$

When $H$ is a polynomial  (in $u$, $\bar u$ and $e^{ix_k}$, $x_k\in[0, 2\pi)$, $k=1, 2, ..., b$) and under the above three assumptions,  the main result is 
\proclaim
{Theorem}
Assume $a=\{a_k\}_{k=1}^b\in (0,\delta]^b=\Cal B(0, \delta):=\Cal B$.  For all $0<\epsilon<1$, there exist $N_0$, $\delta_0>0$ such that 
if $N>N_0$, then for all $\delta\in(0, \delta_0)$, there is  
a Cantor set 
$\Cal G\subset\Cal B$ with 
$$\frac {\text{meas }\Cal G}{\text{meas }\Cal B}\geq 1-\epsilon.$$ 
For all $a\in\Cal G$, there is a quasi-periodic solution of $b$ frequencies to the nonlinear Schr\"odinger equation (1.1):
$$u(t, x)=\sum_{k=1}^b a_k e^{-i{\omega_k}t}e^{ij_k\cdot x}+o(\delta^{2}),$$
with basic frequencies $\omega=\omega(a)=\{\omega_k(a)\}^b_{k=1}$ satisfying 
$$\omega_k=j_k^2+\Cal O(\delta^{2p}),$$
and the amplitude-frequency map $a\mapsto\omega (a)$ is a diffeomorphism. 
The remainder $o(\delta^{2})$ is in a Gevrey norm on $\Bbb T^{b+d}$.
\endproclaim

To our knowledge, this Theorem represents the only known existence result of quasi-periodic solutions 
to the NLS in (1.1) in arbitrary dimension $d$ and for arbitrary nonlinearity $p$.  
When $d>2$, for $p$ sufficiently ``large'', global smooth solutions to (1.1) are not known in general.
These NLS are usually called energy supercritical.   
(Note, moreover, that the equation breaks translation invariance, since $H$ has explicit $x$-dependence.) 

Further, it is seen as a 
Fourier restriction type of theorem on the {\it space-time} phase space. The {\it a priori} lift in (1.4) 
renders much needed flexibility to the method -- in particular -- it is essentially {\it independent} of the
specifics of the underlying linear flow. Indeed, it is applicable to the nonlinear wave equations (NLW) \cite{W4, 5}, which has 
{\it dense} geodesic flow. 

\noindent{\it Remark.}  The bi-characteristics for NLW are hyperbolic, being hyperboloids (or cones). As far as stability
issues are concerned, this presents a more difficult geometry relative to the paraboloids, which are the bi-characteristics for NLS and are limit-elliptic. 
The additional ingredient in \cite{W4, 5} is a Diophantine property of algebraic numbers.

\subheading  {1.2 About the Theorem}

The proof consists of a bifurcation analysis, and an actual construction of the solutions using
a Newton scheme adapted from the work of Bourgain in \cite{B3}, to be explained in more 
details later in this section. The key new aspect is the bifurcation analysis -- to prove the invertibility
of appropriate linearized operators. The NLS equation in (1.1) is completely resonant. When seen as a functional 
equation in the phase space, the linearized operator at $0$ has an infinite dimensional kernel -- the bi-characteristics $\Cal C$.
The idea is then to linearize at the unperturbed solution $u^{(0)}$ instead. The genericity conditions (Gi-iii) in sect. 2 are used to show that 
to leading order, this new linearized operator is a block-diagonal matrix. The determinants of these block matrices are 
polynomials in the amplitude $a$ and control the inverse of the linearized operator at $u^{(0)}$. Here the reasoning is {\it non-perturbative}. 

In general, without additional conditions, one cannot expect an arbitrary linear solution to bifurcate to 
a (nearby) nonlinear solution. The genericity conditions (Gi-iv) in sect. 2 are sufficient conditions for 
such a bifurcation. Technically they bound the sizes of the diagonal blocks of the linearized operator.
Indeed when the genericity 
conditions are violated, qualitatively different dynamics can be exhibited, at least for finite time, 
as  in \cite{CKSTT}, cf. \cite {GuK}.  

The technical conditions ($\flat, \flat\flat$) are only used 
in sect. 4, Lemma 4.2, to prove that these polynomials in $a$ are not identically zero by showing that 
at $a= (1, 1, .., 1)$, the block matrices are diagonally dominant. Hence varying $a$ leads to invertibility
and one can continue with the Newton construction in sect. 5. 

It cannot be excluded that by a more involved combinatorial analysis, particularly for fixed $p$, e.g., $p=1$, these 
two conditions could be improved.  Indeed for the cubic NLS,
actually no restriction on $b$ is needed, therefore the 
Theorem holds for all $b$, cf. the end of the Proof of Lemma 4.2 and the Remark afterwards. 

The loss of $\epsilon$ in the measure estimates is because of the requirement 
$$\Vert (\frac{\partial \omega} {\partial a} )^{-1}\Vert \lesssim\Cal O_{\epsilon}(1) \delta^{2p},$$
to show invertibility. For the cubic NLS, $p=1$, the above bound holds with a {\it uniform} constant $\Cal O(1)$
due to the special polynomial structure of the frequency modulation, cf. (1.18).
So,  in fact,  $\text{meas } \Cal G\to \text{meas } \Cal B$ as $\delta\to 0$. 
Note that in the non-resonant case, the corresponding bound of the parameter-frequency 
modulation is of order $\Cal O(1)$. So the above problem does not arise either and as $\delta\to 0$, 
$\text{meas } \Cal G\to \text{meas } \Cal B$, cf. \cite{B1, 3}.

The further restriction to polynomial nonlinearity $H$ is in order to use directly the analysis in \cite {B3} in sect. 5. 
This is of a technical nature and the Theorem most
likely holds in the analytic category, cf. for example, \cite{sect. 14, B1}. In sects. 2-4, this restriction is not needed. 

The theorem also holds when there is in addition an overall phase, $m\neq 0$, corresponding to adding $mu$ to the right side of (1.1).  
This is because what matters for the bifurcation analysis is the set of differences of the eigenvalues and not the eigenvalues themselves. 

When $d=p=1$, the non-generic set $\Upsilon=\emptyset$. {\it All} $u^{(0)}$ are generic and only amplitude selection is necessary. This is the well understood scenario after writing (1.1) as an infinite dimensional Hamiltonian equation \cite{KP}. 
In this case, the equation is completely integrable.

\subheading {1.3 The cubic NLS}

Among the families of NLS equations in (1.1), the cubic NLS, corresponding to taking $p=1$, has some 
special properties (as we have seen earlier). Quasi-periodic solutions were constructed previously using partial Birkhoff normal forms in dimensions 1 and 2 \cite{B2, GXY, KP}.
This is recently generalized by Procesi and Procesi in \cite{PP1, 2} to arbitrary dimensions for the translationally invariant cubic NLS (the corresponding $H$ has no explicit $x$-dependence). Their result includes linear stability. 
Some of the genericity conditions in \cite{PP1} seem to bear a certain resemblance to the conditions in sect. 2. But the 
formulations of the problem and points of view are rather disparate. In particular, the construction in this paper 
is a general construction in the entire lifted space. The idea behind the subtle definition of the generic linear solutions 
permits our method to bypass the limitations of periodic flow.
\smallskip
We now motivate and describe the method. 

\subheading {1.4 Bifurcation analysis and Lyapunov-Schmidt decomposition}

We express (1.1) using the ansatz in (1.4). By analogy with the standard Fourier series we write $\hat u$ for $\frak a$ and define $\hat {\bar u}$ to be
$\hat {\bar u}(n, j)=\bar \frak a(-n, -j)$ for all $(n,j)\in\Bbb Z^{b+d}$. To simplify notations we write $\hat v$ for $\hat {\bar u}$.
Equation (1.1) can then be written as a nonlinear (infinite) matrix equation: 
$$\text{diag }(n\cdot\omega+j^2)\hat u+(\hat u*\hat v)^{*p}* \hat u+ \sum_{m=1}^\infty\hat\alpha_m*(\hat u*\hat v)^{*(p+m)}* \hat u=0,\tag 1.9$$
where diag $\cdot$  denotes a diagonal matrix, $\omega\in\Bbb R^b$ is to be determined and
$$|\hat\alpha_m(\ell)|\leq C'e^{-c'|\ell|}\quad (C', c'>0)$$
for all $m$. 

From now on we work with (1.9), for simplicity we drop the hat and write $u$ for $\hat u$ and $v$ for $\hat v$ etc.
We seek solutions close to the linear solution $u^{(0)}$ of $b$ frequencies, 
$\text{supp } {u}^{(0)}=\{(-e_{k}, j_k), k=1,...,b\},$ with frequencies
$\omega^{(0)}=\{j_k^2\}_{k=1}^{b}$ ($j_k\neq 0$) 
and small amplitudes $a=\{a_k\}_{k=1}^b$ satisfying $\Vert a\Vert=\Cal O(\delta)\ll 1$.

We complete (1.9) by writing the equation for the complex conjugate. So we have 
$$
\cases
\text{diag }(n\cdot\omega+j^2)u+(u*v)^{*p}* u+\sum_{m=1}^\infty\alpha_m*(u*v)^{*(p+m)}*u=0,\\
\text{diag }(-n\cdot\omega+j^2)v+(u*v)^{*p}* v+\sum_{m=1}^\infty\alpha_m*(u*v)^{*(p+m)}*v=0,
\endcases\tag 1.10
$$
By ``supp", we will always mean the Fourier support, so we write $\text{supp } u^{(0)}$ for 
$\text{supp } {\hat u}^{(0)}$ etc. Let 
$$\aligned\Cal S=&\text{supp } u^{(0)}\oplus\text{supp } {\bar u}^{(0)}\\
:=&S\oplus\bar S.\endaligned\tag 1.11$$

Denote the left side of (1.10) by $F(u, v)$. We make a Lyapunov-Schmidt decomposition into the $P$-equations:
$$ F(u, v)|_{\Bbb Z^{b+d}\times\Bbb Z_2\backslash\Cal S}=F(u, v)|_{\Cal S^c}=0,\tag 1.12$$
and the $Q$-equations:
$$ F(u, v)|_{\Cal S}=0.\tag 1.13 $$
We seek solutions such that 
$u|_S=u^{(0)}$. 
The $P$-equations are infinite dimensional and determine $u$ in the complement of $\text{supp }u^{(0)}$; 
the $Q$-equations are $2b$ dimensional and determine the frequency $\omega=\{\omega_k\}_{k=1}^b$. 

We remark that the above Lyapunov-Schmidt decomposition is the {\it same} as for non-resonant equations, 
in other words, ``as if'' (1.10) were non-resonant. It is a complete change from the usual division of labor, which puts all
the resonances in the $Q$-equations -- the bifurcation equations.
With the decomposition in (1.12, 1.13), the burden of the resonances -- the ``zero-divisors'', 
is shifted to the $P$-equations. Subsequently the resonance analysis -- the bifurcation analysis, is done algebraically, giving rise 
to the genericity conditions in sect. 2. The  ``zero-divisors'' are then turned into small-divisors in sects. 3 and 4,
and (1.10) is transformed into a non-resonant system. Below we give a sketch of this transformation.

We use a Newton scheme to solve the $P$-equations, and as explained earlier, because of the resonances, we linearize the operator
at  $u^{(0)}$,  $v^{(0)}$, instead of at $0$. Let 
$F'( u^{(0)}, v^{(0)})$ be the linearized operator
on $\ell^2(\Bbb Z^{b+d})\times \ell^2(\Bbb Z^{b+d})$,
$$F'=D+A,\tag 1.14$$
where 
$$
D =\pmatrix \text {diag }(n\cdot\omega+j^2)&0\\ 0& \text {diag }(-n\cdot\omega+j^2)\endpmatrix$$
and
$$\aligned
A&=\pmatrix (p+1)(u*v)^{*p}*& p(u*v)^{*p-1}*u*u*\\ p(u*v)^{*p-1}*v*v*& (p+1)(u*v)^{*p}*\endpmatrix+\Cal O(\delta^{2p+2})
\quad  (p\geq 1),\\
&: =\delta ^{2p}A_0+\Cal O(\delta^{2p+2}),\endaligned\tag 1.15$$
with $\omega=\omega^{(0)}$, $u=u^{(0)}$, $v=v^{(0)}$ and we used $\Vert  u\Vert=\Cal O( \delta)$ and homogeneity
to extract the $\delta ^{2p}$ factor in front of $A_0$.

First recall the formal Newton scheme: the first correction
$$\Delta \pmatrix u^{(1)}\\v^{(1)}\endpmatrix= \pmatrix u^{(1)}\\v^{(1)}\endpmatrix- \pmatrix u^{(0)}\\v^{(0)}\endpmatrix
=-[F_{\Cal S^c}'( u^{(0)}, v^{(0)}]^{-1} F_{\Cal S^c}( u^{(0)}, v^{(0)}),\tag 1.16$$
where $F_{\Cal S^c}'$ is $F'$ restricted to $\Cal S^c$: 
$F_{\Cal S^c}'(x, y)= F'(x,y)$, if $x, y\in \Cal S^c$, $F_{\Cal S^c}'(x, y)=0$ otherwise;
$F_{\Cal S^c}(x)=F(x)$ if $x\in \Cal S^c$ and $0$ otherwise.

The operator $A$ is a convolution matrix. Let $P_\pm$ be the projections onto $\Cal C^{\pm}$ and 
$$P=\pmatrix P_+&0\\0&P_-\endpmatrix.\tag 1.17$$ 
Since we look at small data, $\Vert A\Vert=\Cal O(\delta^{2p})\ll 1$ and the diagonals: $\pm n\cdot\omega+j^2$
are integer valued, using ideas from the Schur complement reduction \cite{S1, 2}, the analysis of the spectrum of $F'$ around $0$ can be reduced to that of the projected operator 
$PF'P$ on 
$$\ell^2(\Cal C)=\ell^2(\Cal C^+)\times \ell^2(\Cal C^-),$$
where $\Cal C^+$ and $\Cal C^-$ are as defined in (1.5, 1.6), and to $\Cal O(\delta^{2p+2})$ it is 
the same as the spectrum of $\delta^{2p}PA_0P$ on $\ell^2(\Cal C)$. Thus the heart of the matter is to show that $0\notin\sigma (PA_0P)$ 
on $\ell^2(\Cal C)$.  

We accomplish that by first characterizing the geometry of $u^{(0)}$ in sect. 2.  and then in sect. 3, we show that under 
the genericity conditions (Gi-iii),  $PA_0P$ is a block-diagonal matrix -- a direct sum of matrices in the amplitude $a$ of sizes at most $(2b+2d)\times (2b+2d)$.  
Since $A_0$ is a convolution matrix, $PA_0P$ contains only {\it finite} types of block matrices.
Varying $a$ thus leads to invertibility of each block matrix and therefore $PA_0P$; and consequently $F'_{\Cal S^c}$ for small $\delta$, 
and moreover exponential off-diagonal decay of $[F'_{\Cal S^c}]^{-1}$ ,  using the block structure. Substitute this into (1.16) solves the $P$-equations in the first iteration.

Solving the $Q$-equations in (1.13) gives that the new frequencies $\omega^{(1)}=\{ \omega_k^{(1)}\}_{k=1}^b$ are
$$\omega_k^{(1)}=j_k^2+\frac{{(u^{(0)}*v^{(0)}})^{*p}*u^{(0)}}{a_k}(-e_k, j_k)+\Cal O(\delta^{2p+2}),\, k=1, 2, ..., b.\tag 1.18$$
We note that the $Q$-equations are solved exactly, $F_\Cal S=0$ always. 
Since  $\omega^{(1)}=\omega^{(1)}(a)$, varying $a$ leads to Diophantine frequencies.

From (1.18), the frequency modulation $\Delta\omega^{(1)}$ is of order $\Cal O(\delta^{2p})$, the same order as the matrix $A$
defined in (1.15). So the $P$-equations are still non-amenable. In sect. 4, we iterate the Newton scheme once again
using the modulated $\omega^{(1)}$ as well as the conditions ($\flat$, $\flat\flat$) , and turn the $P$-equations into a non-resonant system. 
In sect. 5,  using the Diophantine frequency $\omega^{(1)}$ and applying the non-resonant analysis scheme in \cite{Chap 19, B3} 
proves the Theorem.  

We note that the Lyapunov-Schmidt method in the present context was introduced in \cite{CW}. It was inspired by the multi-scale analysis in  \cite{FS}.  
The theory was greatly developed by Bourgain to construct quasi-periodic solutions in arbitrary dimensions \cite{B1, 3} and has broad applications. 
More recently, Eliasson and Kuksin \cite{EK} developed a KAM theory in the Schr\"odinger context, which proves existence and linear stability.

All the above results are, however, for 
non-resonant systems, typically using Fourier multipliers as external parameters. In that case, the bifurcation analysis in the first two steps
of the Newton iteration described above can be avoided, as the linearized operator at $0$ is invertible after direct excisions in the parameter space. 

To conclude the introduction, we mention that there is in addition a linear component to the bifurcation theory here. It concerns $L^p$ estimates of $L^2$ eigenfunctions of the Schr\"odinger operator \cite{W3}, 
cf. also \cite{W2}.

\smallskip
\noindent {\it Notations}

We summarize below some of the notational conventions:

\noindent-- The dimension $d$, the degree of nonlinearity $p$ and the number of basic frequencies $b$ are fixed. The set 
$\{j_k\}_{k=1}^b$ is a fixed subset of $\Bbb Z^d$.

\noindent-- The letter $u$ denotes a function on $\Bbb T^{b+d}$, $\hat u$ its Fourier series and $\hat v$ the Fourier 
series of $\bar u$.  One generally drops the hat and writes $u$ for $\hat u$ and $v$ for $\hat v$, which are 
functions on $\Bbb Z^{b+d}$. 

\noindent-- The letters $n$ and $\nu$ denote vectors in $\Bbb Z^b$; while $j$ and $\eta$ vectors in $\Bbb Z^d$.  

\noindent-- The dot $\cdot$ denotes the usual scalar product in Euclidean space. To simplify notations, one writes $j^2$ for $j\cdot j$ etc. 

\noindent-- The norm $\Vert\,\Vert$ stands for the $\ell^2$ or operator norm; while $|\,|$ for the sup-norm or the length of a vector in a finite 
dimensional vector space or the number of elements in a given set.

\noindent-- Given two sets $A$ and $B$, the sum $A+B$ is defined in the usual way. If $A=B$, then one
writes $A+A=2A$; consequently $2A+A=3A$ etc.

\noindent-- A matrix of vectors is denoted by $[[\quad ]]$.

\noindent-- An identically zero function $f$ is denoted by $f \equiv  0$; the negation $f\not\equiv 0$.
\smallskip
\noindent{\it Acknowledgement.} The author wishes to express her gratitude to the referees, whose extraordinary
care and thoroughness, in reading the manuscript, greatly helped the exposition of the main ideas. This work was partially
supported by the grant ANR-10-JCJC0109.
\bigskip 
\head{\bf 2. The genericity conditions}\endhead
In this section,
genericity conditions will be imposed on the spatial Fourier support $\{j_k\}_{k=1}^b$ of $u^{(0)}$. 
This paves the way toward showing in sect. 3 that 
$PA_0P$ defined by (1.17, 1.15) reduces 
to a block diagonal matrix with finite types of blocks. As will be amplified later in the section, 
the genericity conditions stem from bounding the sizes of these block matrices.
The blocks are described
algebraically by linear and quadratic 
polynomials in the spatial Fourier variable $j\in\Bbb Z^d$ with {\it coefficients} dependent on 
$\{j_k\}_{k=1}^b$ -- the size of a block is bounded above by the number 
of possibly compatible equations.  The (to be stated) genericity conditions will
yield precise upper bounds on the sizes of such blocks in sect. 3.  

\subheading {2.1 Basic notions}

We make the identification:  
$$\Bbb Z^{b+d}\times\Bbb Z_2\sim\{ \Bbb Z^{b+d}, +\}\cup \{\Bbb Z^{b+d}, -\}.\tag 2.1$$
We write a matrix $M$ acting on $\Bbb Z^{b+d}\times\Bbb Z_2$:
$$\pmatrix M_{11}&M_{12}\\{M_{21}}&{M_{22}}\endpmatrix,$$
where $M_{ij}$, $i$, $j=1$, $2$, act on $\Bbb Z^{b+d}$, as
$$M=M_{11}\otimes\pmatrix 1&0\\0&0\endpmatrix + M_{12}\otimes\pmatrix 0&1\\0&0\endpmatrix+M_{21}\otimes\pmatrix 0&0\\1&0\endpmatrix
+M_{22}\otimes\pmatrix 0&0\\0&1\endpmatrix.\tag 2.2$$

To define generic $u^{(0)}$, we need to analyze the convolution
matrix $A_0$ in (1.15). We use the notation introduced in (1.3). 
$u^{(0)}$ has support: 
$$\text { supp } u^{(0)}=\{(-e_k, j_k), k=1, 2, ..., b\},$$
where $e_k\in\Bbb Z^b$ are the basis vectors, $j_k\in\Bbb Z^d$ are fixed vectors and 
$(-e_k, j_k)\in\Bbb Z^{b+d}$.

Let $S$ be the set of {\it ordered} pairs:
$$S=\{(k, k'), k\neq k', k, k'=1, 2, ..., b\}.$$
For $s=(k, k')\in S$, define the difference of pairs:  
$$\aligned \nu_s&=e_k-e_{k'}\in\Bbb Z^b,\\
\eta_s&=j_k-j_{k'}\in\Bbb Z^d.\endaligned\tag 2.3$$
Note that  the map $s\to \nu_s$ is {\it invertible} on the range of $\nu_s$.
Using (2.3), we then have 
$$\text { supp } |u^{(0)}|^2= \{(-\nu_s, \eta_s), s\in S\}\bigcup \{(0, 0)\}\subset\Bbb Z^{b+d}.\tag 2.4$$

Fix $p\in\Bbb N^+$. Denote by $\{p_s\}_{s\in S}$, a family in $\Bbb N$ satisfying 
$$p_s\geq 0 \text { and }\sum_{s\in S} p_s\leq p.\tag 2.5$$
For each family $\{p_s\}_{s\in S}$, define 
$$\align \nu&=-\sum_{s\in S} p_s\nu_s\in\Bbb Z^b,\tag 2.6\\
 \eta&=\sum_{s\in S} p_s\eta_s\in\Bbb Z^d.\tag 2.7\endalign$$
Denote by $\{(\nu, \eta)\}$ the set of all $(\nu, \eta)$ constructed by (2.5-2.7).
Fourier series multiplication then gives
$$\Gamma:=\text{ supp } |u^{(0)}|^{2p}=\{(\nu, \eta)\}\subset\Bbb Z^{b+d}.\tag 2.8$$
This takes care of the combinatorial description of the symbol of convolution matrices $(A_0)_{11}$ and $(A_0)_{22}$
in (1.15).

Further, for later constructions, it is convenient 
to view $\eta$ as maps. One first makes the identification
$\{j_k\}_{k=1}^b\sim (j_1, j_2, ..., j_b)\in(\Bbb Z^d)^b$.
From (2.3),  for all $\nu_s$, $s\in S$, there is the relation:
$$\eta_s=\sum_{i=1}^b\nu_s^ij_i,    \tag $\sharp$ $$
where $\nu_s^i$ is the $i$th component of $\nu_s$. So using the 
invertibility of the map $s\to\nu_s$, 
for each given $s$, $\eta_s$ 
defines a map from $(\Bbb Z^d)^b\to\Bbb Z^d$:   
$$\eta_s:\, (\Bbb Z^d)^b\ni(j_1, j_2, ..., j_b)\mapsto j_k-j_{k'}\in\Bbb Z^d.$$
Therefore for a given family $\{p_s\}$ satisfying (2.5), 
the $\eta=\eta_{\{p_s\}}$ constructed in (2.7) defines 
a map from $(\Bbb Z^d)^b$ to $\Bbb Z^d$. The $\nu=\nu_{\{p_s\}}$ in (2.6)
are constants as the $\nu_s$.

For example, take $d=1$, $b=3$ and write $(x, y, z)$ for a vector in $\Bbb Z^3$: 
$(x, y, z)\in\Bbb Z^3$. If $s=(1,2)$, then $$\nu_s=(1, -1, 0)\in\Bbb Z^3,$$ and 
$$\eta_s(x, y, z)=x-y\in\Bbb Z$$ is a map from $\Bbb Z^3$ to $\Bbb Z$. 

\proclaim{Lemma 2.1} For a given family $\{p_s\}$, $$(\nu, \eta)=(\nu_{\{p_s\}}, \eta_{\{p_s\}})\in \text{supp }|u^{(0)}|^{2p}$$
satisfies the following relations:
$$\nu=0 \Longrightarrow \eta\equiv 0$$ and 
$$\eta\equiv 0\Longrightarrow \nu=0.$$
So
$$\nu=0 \Longleftrightarrow\,  \eta\equiv 0.\tag 2.9$$
\endproclaim
\demo{Proof}  
From (2.6, 2.7), there is the analogue (generalization)  of $(\sharp)$: 
$$\eta=-\sum_{i=1}^b\nu^ij_i,    \tag $\sharp\sharp$ $$
where $\nu^i$ is the $i$th component of $\nu$. 
Clearly if $\nu=0$, then $\eta\equiv 0$. If $\eta\equiv 0$ and $\nu\neq 0$, then there exists 
$k\in\{1, 2, ...b\}$ such that $\nu^{k}\neq 0$. Set all $j_{k'}=0$ for $k'\neq k$, then 
$\eta\neq 0$ for $j_k\neq 0$, which is a contradiction.  \hfill $\square$

\enddemo
\noindent{\it Remark.} The constant vectors $\nu=\nu_{\{p_s\}}$ do the ``book keeping"
when multiplying Fourier series.  

We also need to analyze the Fourier support of the other two symbols $(A_0)_{12}$ and $(A_0)_{21}$
in (1.15). They are complex conjugates of each other. We have similarly
$$\Gamma_+:= \text{ supp } |u^{(0)}|^{2(p-1) }[u^{(0)}]^2=\{(\nu, \eta)\}\subset\Bbb Z^{b+d}\tag 2.10$$
where 
$$\aligned \nu&=-\sum_{s\in S} p_s\nu_s-(e_\kappa+e_{\kappa'}),\\
 \eta&=\sum_{s\in S} p_s\eta_s+(j_\kappa+j_{\kappa'}),\endaligned\tag 2.11$$
with $$\sum_{s\in S} p_s\leq p-1,\, p_s\geq 0, \kappa, \kappa'=1, 2, ..., b;$$
and 
$$\Gamma_-:= \text{ supp } |u^{(0)}|^{2(p-1) }[v^{(0)}]^2=\{(\nu, \eta)\}\subset\Bbb Z^{b+d}\tag 2.12$$
where 
$$\aligned \nu&=-\sum_{s\in S} p_s\nu_s+(e_\kappa+e_{\kappa'}),\\
 \eta&=\sum_{s\in S} p_s\eta_s-(j_\kappa+j_{\kappa'}),\endaligned\tag 2.13$$
with $$\sum_{s\in S} p_s\leq p-1,\, p_s\geq 0, \kappa, \kappa'=1, 2, ..., b.$$
We note that $\nu\neq 0$ and $\eta\not\equiv 0$ in (2.11) and (2.13).

For the definition below, it is again more convenient to view $(\nu, \eta)$ as a point in $\Bbb Z^{b+d}$.
Define the sets $W_r\subset\Bbb Z^{b+d}$, $r=1, 2, ...$ as follows:
$$\aligned W_1&=\Gamma\bigcup \{\Gamma+\Gamma_-\},\\
W_2&=2\Gamma\bigcup \{2\Gamma+\Gamma_-\},\endaligned$$
$$\vdots$$
$$W_r=r\Gamma\bigcup \{r\Gamma+\Gamma_-\}\tag 2.14$$
$$\vdots$$
We note that 
$W_r\subset W_{r'},$ if $r<r'$ and 
since $(0, 0)\notin\Gamma_-$, $$r\Gamma\bigcap \{r'\Gamma+\Gamma_-\}=\emptyset\tag 2.15$$
for all $r$, $r'=1, 2, ...$ Note also that for all $r$, all elements $(\nu, \eta)$ in $W_r$ satisfy the relation in ($\sharp\sharp$).

\noindent{\it Remark.} As will become clear later in the section, the reason for the definition of $W_r$ is that both sets in the union generate
$r$ possible linear equations in $j\in\Bbb Z^d$.    

To analyze further these sets, we define
$$\Cal W_R=W_R\backslash\{(0, 0)\}, \, R=1, 2, ...\tag 2.16$$
and 
$$\Cal W'_R=R\Gamma\backslash\{(0, 0)\}\subset\Cal W_R, R=1, 2, ...\tag 2.17$$
Elements of $\Cal W_R$ are again denoted by $(\nu,\eta)$. 

To gain a better insight into these sets, below are some examples of $(\nu, \eta)$ in $\Cal W'_1$, $\Cal W'_2$, $\Cal W'_3$ and $\Gamma_-$ 
in the case $p=1$. For simplicity, we only display the $\eta$ component. 

\noindent $\Cal W'_1$: $j_k-j_{k'}$, $k\neq k'$, $k, k'=1, 2, ..., b$;

\noindent $\Cal W'_2$:  $2(j_1-j_2)$, .., $j_1-2j_2+j_3$, ..., $j_1-j_2+j_3-j_4$, ... ;

\noindent $\Cal W'_3$:  $3(j_1-j_2)$, ..., $2(j_1-j_2)+j_3-j_4$, ..., $j_1-2j_2+j_3+j_5-j_6$, ..., $j_1-j_2+j_3-j_4+j_5-j_6$, ... ;
  
\noindent $\Gamma_-$: $-(j_k+j_{k'}), k, k'=1, 2, ..., b$.  

We resume the construction and fix an $R$ and let $\sigma$ be a subset of $\Cal W_R$:  
$$\sigma\subset \Cal W_R\subset\Bbb Z^{b+d}.$$
Denote by $|\sigma|$ the number of elements in $\sigma$ and call $|\sigma|$ the length. 
We first define a notion of a {\it connected} subset of $\Cal W_R$. 

Given a subset $\sigma$ of $\Cal W_R$, in view of (2.15) and since the problem is posed on $\Bbb Z^{b+d}\times\Bbb Z_2$, 
we define a map 
$$\aligned \sigma&\mapsto \Bbb Z^{b+d}\times \Bbb Z_2:\\
\sigma\ni\varsigma&\mapsto \varsigma'\in \Bbb Z^{b+d}\times \Bbb Z_2;\endaligned$$
for 
$$\varsigma\in\{\sigma\cap\Cal W'_R\}\subset  \Cal W_R\subset\Bbb Z^{b+d};$$
 define 
 $$\varsigma':=(\varsigma,+)\in \Bbb Z^{b+d}\times \Bbb Z_2;\tag 2.18$$
 for 
 $$\varsigma\in\sigma\cap\{\Cal W_R\backslash\Cal W'_R\}$$
 define
 $$\varsigma':=(\varsigma,-).\tag 2.19$$
 Let $$\sigma':=\{\varsigma'\}\cup \{(0, 0,+)\},\tag 2.20$$
where $\{\varsigma'\}$ denotes the set formed from $\sigma$ as above.

\noindent{\bf Definition.}
A set $\sigma\subset \Cal W_r\subset \Bbb Z^{b+d}$ is {\it connected} if the image set  $$\sigma'=\sigma'(\sigma)\subset \Bbb Z^{b+d}\times \Bbb Z_2$$ has the
following properties and is therefore also called {\it connected}: 

\noindent For all $a'$, $b'\in\sigma'$, $\exists \varsigma'_1, \varsigma'_2, ..., \varsigma'_k\in\sigma'$,
 such that if we set $\varsigma'_0=a'$ and $\varsigma'_{k+1}=b'$, then 
 if $$\varsigma'_{i+1}=(\varsigma_{i+1},+)\text{ and }\varsigma'_{i}=(\varsigma_{i},+)\tag 2.21$$
 or if $$\varsigma'_{i+1}=(\varsigma_{i+1},-)\text{ and }\varsigma'_{i}=(\varsigma_{i},-),$$
 then $$\varsigma_{i+1}-\varsigma_i\in\Gamma\backslash\{(0, 0)\};$$
 if  $$\varsigma'_{i+1}=(\varsigma_{i+1},+)\text{ and }\varsigma'_{i}=(\varsigma_{i},-),$$
 then $$\varsigma_{i+1}-\varsigma_i\in\Gamma_{+};$$
 or else 
 if  $$\varsigma'_{i+1}=(\varsigma_{i+1},-)\text{ and }\varsigma'_{i}=(\varsigma_{i},+),$$
then  $$\varsigma_{i+1}-\varsigma_i\in\Gamma_{-}.\tag 2.22$$
 
 \noindent{\it Remark.} If $\sigma\subset \Cal W'_R\subset\Bbb Z^{b+d}$, one could simplify and define 
 $\sigma'=\{(0, 0)\}\bigcup \sigma\subset\Bbb Z^{b+d}$. One then has that $\sigma$ is a connected set 
 if $$\forall a, b \in\sigma',\exists \varsigma_1, \varsigma_2, ... \varsigma_k\in\sigma',$$
 such that if we set $\varsigma_0=a$ and $\varsigma_{k+1}=b$, then $\varsigma_{k+1}-\varsigma_k\in\Gamma\backslash\{(0, 0)\}$.
 
Following directly from the definition, there is 
\proclaim {Lemma 2.2} For any connected set of length at least $2$: $\sigma\subset\Cal W_R\subset\Bbb Z^{b+d}$ with $|\sigma|\geq 2$ and $\sigma\cap\Cal W'_R\neq\emptyset$, 
define $\tilde\sigma=\sigma\cap\Cal W'_R$ and $\bar\sigma=\{(0,0)\}\cup\tilde\sigma$. Then 
$$\forall a, b \in\bar\sigma,\exists \varsigma_1, \varsigma_2, ... \varsigma_k\in\bar\sigma,$$
 such that if we set $\varsigma_0=a$ and $\varsigma_{k+1}=b$, then $\varsigma_{\ell}-\varsigma_{\ell'}\in |\sigma|\Gamma\backslash\{(0, 0)\},$
 for all $\ell, \ell'=0, 1, 2, ..., k+1$, $\ell\neq \ell'$.
\endproclaim
\demo{Proof} Clearly $$\varsigma_{\ell}-\varsigma_{\ell'}\in r \Gamma\backslash\{(0, 0)\}$$
for some $r\in\Bbb N$,  since they are both in $\bar \sigma$. Moreover $r$ must satisfy $r\leq |\sigma|$ 
as $\sigma$ is a {\it connected} set and $$\Gamma+\Gamma=2\Gamma=\Gamma_++\Gamma_-,$$
from (2.8, 2.10, 2.12), cf. also the examples given earlier.
$\square$
\enddemo

 We now fix $R=2(d+1)$ and specialize to 
 $$\Cal W:=\Cal W_{2(d+1)}\subset\Bbb Z^{b+d}.$$ We shall treat the two copies of 
 $\Bbb Z^{b+d}$ in 
 $$\Bbb Z^{b+d}\times \Bbb Z_2\sim \{\Bbb Z^{b+d}, +\}\cup \{\Bbb Z^{b+d}, -\}$$
 ``separately" in some sense. Define $$\Cal W':=\Cal W'_{2(d+1)}=\text{ supp } |u^{(0)}|^{4p(d+1)}\backslash\{(0,0)\},\tag 2.23$$
from (2.17, 2.6, 2.7). 

For each $(\nu,\eta)\in\Cal W'$, define 
$$J=J(\nu,\eta)=|\eta|^2+\nu\cdot\omega^{(0)}:=\eta^2+\nu\cdot\omega^{(0)}\in\Bbb Z,\tag 2.24$$
where $\omega^{(0)}=\{j_k^2\}_{k=1}^b$, as above (1.3), 
and the $d$-dimensional hyperplane in $\Bbb R^d$: 
$$2\eta^1x_1+2\eta^2x_2+...+2\eta^dx_d+J=0,\tag 2.25$$
where $\eta^k$, $k=1, 2,..., d$ is the $k$th component of $\eta\in\Bbb Z^d$ and 
$x_k$, $k=1, 2,..., d$ is the $k$th component of $x\in\Bbb R^d$. For simplicity and by an abuse of 
notation, we shall use the $(d+1)$-dimensional row vector $(2\eta, J)$ to 
denote the above plane. 

\subheading {2.2 The genericity conditions}

Recall from sect. 1 that we start from the linear solution $u^{(0)}$ of $b$ frequencies: 
$$\aligned u^{(0)}(t, x)&=\sum_{k=1}^b a_ke^{-ij_k^2t}e^{ij_k\cdot x}\\
:&=\sum_{k=1}^b\hat u(-e_{k}, j_k)e^{-i(e_{k}\cdot\omega^{(0)}) t}e^{ij_k\cdot x}, \endaligned
$$
where $e_{k}\in\Bbb Z^b$ is the $k$th basis vector, 
$\omega^{(0)}=\{j_k^2\}_{k=1}^{b}$ ($j_k\neq 0$) and $\hat u^{(0)}(-e_{k}, j_k)=a_k$. Below we specify
conditions on the spatial frequencies $\{j_k\}_{k=1}^b\in(\Bbb Z^d)^b$ for $u^{(0)}$ to be generic, in order to construct
nearby nonlinear solutions in sects. 3-5.
\bigskip
\noindent{\bf Definition.} $u^{(0)}$ of $b$ frequencies is {\it generic} if its Fourier support $\{(-e_m, j_m)\}^b_{m=1}\subset\Bbb Z^{b+d}$
satisfies:
\item{(Gi)} For all $j_k$, $k=1, 2, ..., b$, define the set of differences 
$$L_k=\{j_{k'}-j_k|k'=1, ..., b, k'\neq k\}.$$ 
If $b\geq d+1$, any $d$ vectors in $L_k$ are linearly independent.
(If $b\leq d$, there is no condition (Gi).) 
\item{(Gii)} For all $(\nu, \eta)\in\Cal W'$ defined in (2.23), 
$$\eta\neq 0.$$
\item{(Giii)}  For any connected set $\sigma=\{(\nu, \eta)\}\subset\Cal W$, define
$\tilde\sigma=\sigma\cap\Cal W'$. Assume $|\tilde\sigma|=d+1$. (If $|\tilde\sigma|\leq d$, there is no condition (Giii).)
Denote by $\frak{h}:=\{\eta_i\}_{i=1}^{d+1}$, the 
set of $(d+1)$ $\eta$'s in $\tilde\sigma$ and $(2\frak{h}, \Cal J):=\{(2\eta_i, J_i)\}_{i=1}^{d+1}$ the corresponding set of $(d+1)$ planes
in the form (2.25). 
Assume that $\frak h$ is included ($\subseteq$) in {\it no} $L_k$, 
then $$\bigcap _{i=1}^{d+1}(2\eta_i, J_i)=\emptyset.\tag *$$
\item{(Giv)} For all $j_m$, $m=1,..., b$ and all $(\nu, \eta)\in\Gamma\backslash\{(0,0)\}$,  $f\in\Bbb Z$ defined as 
$$f:=\nu\cdot\omega^{(0)}+2j_m\cdot\eta+\eta^2\neq 0,$$
 if $(\nu, \eta)\neq (-e_{k'}+e_m, j_{k'}-j_m)$, for all $k'=1,..., b$.
\bigskip
\noindent{\it Remarks.} 1. The main observation which leads to the formulation of the generic conditions is that 
due to curvature there is a lack of translation invariance on the bi-characteristics. Indeed 
elliptic or limit-elliptic bi-characteristics are specific instances. For the Schr\"odinger equation
here, this is expressed through the functions:
$\eta=\eta (j_1, j_2, ..., j_b)$ defined in  (2.7) and $J=J(j_1,  j_2, ..., j_b)$ defined in (2.24).
 We have that $\eta$ is {\it linear} in  $j_1$,  $j_2$, ...,~ $j_b$; while $J$ 
{\it quadratic}.  

\noindent 2. To understand the conditions (Gi-iii), we preview here that there is a corresponding notion
of connected sets on the bi-characteristics $\Cal C$, to be defined in ($\dag$) later 
in the section. These connected sets give rise to the block matrices in $PA_0P$, mentioned at
the beginning of the section.  Conditions (Gi, ii) ensure that the exceptional connected sets on $\Cal C$ are of size $2b$ and 
have spatial support the set $\{\pm j_k\}_{k=1}^b$. Moreover under the additional condition on the frequencies in ($\flat$), they are all sufficiently far 
from the origin and lie outside the domain of the $P$-equations for the second Newton iteration in sect. 4. Condition (Gii) ensures that the 
``hyperplanes'' appearing in (*) are proper subsets of co-dimension $1$.  Condition (Giii) then isolates the exceptional sets and ensures that all other connected sets have sizes at most 
$(2d+2)$. These bounds will be proven in Lemma 3.2, sect. 3 -- they are an essential ingredient for the dynamics exhibited in the Theorem.
``Large'' connected sets could possibly lead to qualitatively different behavior, as mentioned in sect. 1. 

\noindent 3. Condition (Giv) indicates that for all $m=1,2 , ..., b$,  the frequency vectors $j_m$
do not lie in the planes defined by $(0,0) \neq (\nu, \eta)\in \Gamma\backslash L_m$.
This ensures that the Lyapunov-Schmidt 
$P$ and $Q$ equations decomposition is stable for the first two iterations (and consequently all subsequent iterations).
In other words, the ``normal'' directions are normal to the ``tangential'' directions for all time.
From reflection symmetry, it suffices to consider the $f$ defined there. 
\smallskip
Since $J$, $\eta$, $f$ are viewed as functions of $\{j_k\}_{k=1}^b$ and so will be the various determinants used to describe (*) in condition (Giii),  
the genericity conditions above are algebraic conditions on $\{j_k\}_{k=1}^b$.
Let $\Upsilon_1$ be the set on which the first genericity condition (Gi) is violated and similarly define 
$\Upsilon_k$ for $k=2, 3, 4$.
The following lemma fulfills the prerequisite for the construction.
\proclaim{Lemma 2.3}
The non-generic set $$(\Bbb R^d)^b\supset\Upsilon:=\bigcup_{k=1}^4\Upsilon_k$$
is algebraic, has co-dimension $1$ and $\{(\Bbb R^{d})^b\backslash \Upsilon\}\cap(\Bbb Z^d)^b$ is an infinite set.
\endproclaim

The proof of Lemma 2.3 will be decomposed into the proofs of 

\proclaim{Lemma 2.4}
The non-generic set $$(\Bbb R^d)^b\supset\Upsilon':=\bigcup_{k\neq 3}\Upsilon_k$$
is algebraic and has co-dimension $1$.
\endproclaim

\proclaim{Lemma 2.5}
The non-generic set $$(\Bbb R^d)^b\supset\Upsilon_3$$
is algebraic and has co-dimension $1$.
\endproclaim
\smallskip
\demo{Proof of Lemma 2.4} To prove the assertions, it suffices that the various algebraic functions defined in (Gi, ii, iv) are non-zero functions, 
which in turn suffices if they are non-zero when restricted to appropriate sub-varieties.

\noindent (i) Since $j_m\in\Bbb R^d$,  $m=1, 2, ..., b$,
$$D=\text{det } [[\{j_{k_i}-j_m\}_{i=1}^d]]\not\equiv 0, \, k_i\neq m,$$
$$D=0$$ 
define sets of co-dimension $1$ in $(\Bbb R^d)^b$.

\noindent (ii): This follows from Lemma 2.1 with the exponent $2p$ replaced by $4p(d+1)$. 

\noindent (iv): (a) If $\eta$ is not a function of $j_m$, restrict to $j_m=0$. 
From (2.6, 2.7), $\nu$ and $\eta$ are functions of at least $2$ variables in the set 
$\{j_\ell\}_{\ell=1}^b$. So  $J=J (j_1, j_2, ..., j_b)$ is a function of at least $2$ variables in the set $\{j_\ell\}_{\ell=1}^b$ and there is $k\in\{1, 2, ..., b\}$ such that when restricting to the sub-manifold in $(\Bbb R^{d})^b$:
$$j_\ell=0, \text{ for }\ell\neq k,$$ 
$\nu\cdot\omega^{(0)}=|\nu\cdot\omega^{(0)}|$ is positive, so there exists $A>0$
$$J(0, ..., j_k, ...,0)=Aj_k^2\not\equiv 0.$$ 

If $\eta$ is a function of $j_m$,  then from the
structure of $\eta$ in (2.7) and the restriction on $\eta$, either 

\item{(b)} $\eta$ is a function of at least $2$ other variables in the set $\{j_k| k=1, ..., b, \, k\neq m\}$

\noindent or 

\item{(c)} it is a function of only  $1$ other variable $j_{m'}$, $m'\neq m$, then the coefficient in front of $j_{m'}$ is not $1$.

Restricting to $j_m=0$, 
\item{(b)} $f$ is a function of at least $2$ variables in the set $\{j_k| k=1, ..., b, \, k\neq m\}$, so this reduces to case (a). 
\item{(c)} there exists $A>0$ such that  $f (0, ..., j_{m'}, ..., 0)=Aj_{m'}^2\not\equiv 0$.

Combining the above, we have proven that $\Upsilon'$ is algebraic of co-dimension $1$. 
$\square$
\enddemo

\demo{Proof of Lemma 2.5}
We give an algebraic description of the geometry entailed by (*) as follows. For a given set $\frak h$, we study the  
sub-matrices of the $(d+1)\times (d+1)$ matrix $[[(2\frak h, \Cal J)]]_{(d+1)\times (d+1)}$.  The goal is to show that
there is a sub-matrix with determinant  $D\not \equiv 0$ and setting $D\neq 0$ yields the 
geometry in (*).  

There are 2 cases: 
\item{a)} all subsets $\frak h'\subset \frak h$ of $d$ vectors, $|\frak h'|=d$, satisfy $\det [[\frak h']]_{d\times d}\not\equiv 0$;
\item {b)} there exists a subset $\frak h'\subset \frak h$, $|\frak h'|=d$, such that $\det [[ \frak h']]_{d\times d}\equiv 0$.

Case a) is further divided into 2 subcases:
\item{a1)} there exists $i\in\{1, 2, ..., d+1\}$, such that $\eta_i$ is a function of at least $3$ variables in $\{j_k\}_{k=1}^b$;
\item{a2)} all $\eta_i$, $i=1, 2, ..., d+1$ are functions of $2$ variables only.

a1) Without loss of generality, one may assume  $i=d+1$, after a possible relabelling. Denote by $j_m$, $m\in\{1, 2, ..., b\}$, a variable that $\eta_{d+1}$ depends on. 
Define $$(2\tilde \eta_i, \tilde J_i)=(2\eta_i, J_i)-c_i(2\eta_{d+1}, J_{d+1}),$$  for $i\neq d+1$, where $c_i\in\Bbb Q$, so that 
$\tilde \eta_i$, $ i=1, 2, ..., d$, are {\it independent} of $j_m$.  
(If $\eta_i$, $i=1, 2, ..., d$, are independent of $j_m$, then $c_i=0$ and $(2\tilde \eta_i, \tilde J_i)= (2\eta_i, J_i)$.)
Define  $$(2\tilde \eta_{d+1}, \tilde J_{d+1})=(2\eta_{d+1},  J_{d+1}).$$
One has  
$$D=\det [[(2\frak h, \Cal J)]]_{(d+1)\times (d+1)}=\det [[(2\tilde \frak h, \tilde \Cal J)]]_{(d+1)\times (d+1)},$$
where $(2\tilde \frak h, \tilde \Cal J)=\{(2\tilde \eta_i, \tilde J_i)\}_{i=1}^{d+1}$. 

There are $2$ possibilities:
\item{(i)}  $$\tilde D=\det [[\{\tilde \eta_i\}_{i=1}^d]]_{d\times d}\not \equiv 0,$$
{\it as a function of} $j_k$, $k\neq m$, 
(recall that by construction $\tilde \eta_i$, $i=1, 2, ..., d$, are only functions of $j_k$, $k\neq m$).

Restrict $D$ to the sub-variety $V$ defined by $\tilde \eta_{d+1}= \eta_{d+1}=0$ and consider the restriction as a function of $j_k$, $k\neq m$. 
The determinant expansion then gives 
$$D=\pm \tilde J_{d+1}\det [[\{2\tilde \eta_i\}_{i=1}^d]]_{d\times d} =\pm 2^dJ_{d+1}\tilde D$$
(By an abuse of notation, we have omitted the restriction sign.) 

Assume that $s$, $3\leq s\leq b$, components of $\nu_{d+1}$ are non-zero. (Here Lemma 2.1 is used to obtain the lower bound $s\geq 3$.) 
Rename these components $\nu^1$, $\nu^2$, ..., $\nu^s$ and rename the $j_i$ accordingly, if necessary, we have on $V$, 
$$\eta_{d+1}=\sum_{i=1}^s \nu^i j_i=0.$$
So $$j_s=-\frac{\sum_{i\neq s}\nu^ij_i}{\nu^s}.$$
Using the above gives 
$$J_{d+1}=\sum_{i=1}^s\nu^i j_i^2=\sum_{i\neq s}\nu^ij_i^2+\frac{(\sum_{i\neq s}\nu^ij_i)^2}{\nu^s}\not\equiv 0,$$
as the last expression contains cross terms such as $j_1\cdot j_2$ etc. 
So $D\not\equiv 0$ on $V$ using also that $\tilde D\not\equiv 0$. So $D\not\equiv 0$. 
On the set defined by $D\neq 0$, (*) is satisfied. 

\item{(ii)} $$\tilde D=\det [[\{\tilde \eta_i\}_{i=1}^d]]_{d\times d}\equiv 0.$$ 
In this case, there must exist $i\in \{1, 2, ..., d\}$, such that $c_i\neq 0$ (as otherwise it contradicts the definition of case a)), and there are  {\it constants} $\alpha_k\neq 0$, $k=1, 2, ..., d+~1$, such that 
$$\sum_ {k=1}^{d+1} \alpha_k\eta_k\equiv 0.$$
Using (2.9), this in turn gives 
$$\Bbb Z^b\ni\sum_{k=1}^{d+1}\alpha_k\nu_k=0.$$
So 
$$\sum_{k=1}^{d+1}\alpha_k\nu_k\cdot\omega^{(0)}\equiv 0$$
and we have 
$$\sum_{k=1}^{d+1}\alpha_kJ_k=\sum_{k=1}^{d+1}\alpha_k|\eta_k|^2.$$
Solving for $\eta_{d+1}$ using the linear relation, we have 
$$\eta_{d+1}\equiv -\sum_{k=1}^d\frac{\alpha_k}{\alpha_{d+1}}\eta_k, \quad \alpha_k\neq 0, \,k=1, 2, ..., d+1.$$
So $$\sum_{k=1}^{d+1}\alpha_k\eta_k^2\equiv \sum_{k\neq d+1}\alpha_k\eta_k^2+\frac{1}{\alpha_{d+1}}(\sum_{k\neq d+1}\alpha_k\eta_k)^2\not\equiv 0,$$ 
as the second sum yields cross terms. 

More precisely, we may set $\alpha_{d+1}=-1$ without loss of generality. Then 
the right side of  the above equation is equal to 
$$
R:=\sum_{k\neq d+1}\alpha_k\eta_k^2-\sum_{k\neq d+1}\alpha_k^2\eta_k^2-2\sum_{m\neq n\neq d+1}\alpha_m\alpha_n\eta_m\cdot\eta_n.
$$
If there exists $\alpha_k\neq 1$, then setting all $\eta_{k'}=0$ for $k'\neq k$.
$$R=(\alpha_k-\alpha_k^2)\eta_k^2\not\equiv 0.$$ 
(Here we also used the determinant conditions defining case a).)
If all $\alpha_k=1$, then setting $2$ of the $\eta_k$ to be $1$ and the rest $0$ give that
$$R=-2\neq 0,$$ which proves the claim. 

So when $$\sum_{k=1}^{d+1}\alpha_k\eta_k\equiv 0,$$
$$\sum_{k=1}^{d+1}\alpha_kJ_k\not\equiv 0.$$
Since $$\det [[\{\eta_i\}_{i=1}^d]]_{d\times d} \not\equiv 0$$
from the definition of case a), 
$$D\not\equiv 0$$ after expanding the determinant. On the set defined by $D\neq 0$, (*) is satisfied. (Note that 
on the same set, $\det [[\{\eta_i\}_{i=1}^d]]_{d\times d} \neq 0$.)

a2) There are 2 sub-cases:
\item{(i)} There exist $i$, $k$, $i\neq k$, such that $\eta_i-\eta_k$ is a function of at least $3$ variables.

Without loss of generality, one may assume  $k=1$ (after a possible relabelling). Subtract the $i$th, all $i\neq 1$, equations of the form (2.25) from the first. The change of variable : $x\to x+\eta_1$ 
subsequently transforms the resulting system of linear equations to case a1) or the to be treated case b). 

\item{(ii)} $\eta_i-\eta_k$ is a function of $2$ variables for all $i\neq k$. 

It follows that there must be $q\in\Bbb Z\backslash \{0, 1\}$, $m\in \{1, 2, ..., b\}$ such that all $\eta_i$, $i=1, 2, ..., d+1$, are of the form 
$$q(j_\ell-j_m), $$ $\ell\neq m$, $\ell =1, 2, ..., b$; $q=1$ is excluded by the condition on $\frak h$ in (Giii).   After relabelling, one 
may assume that $m=1$ and that $\eta_1= q(j_2-j_1)$. Expand the determinant using the first row. Since $(\eta_1, J_1) $ is the 
{\it only} vector that depends on $j_2$, using that $\eta_1$ is linear in $j_2$ and that $J_1$ is {\it quadratic}, 
$D\neq 0$ when $|j_2|\gg 1$ and $\det [[\{\eta_i\}_{i=2}^{d+1}]] \neq 0$. So $$D\not\equiv 0,$$ and on the set defined by $D\neq 0$, (*) is satisfied.
\smallskip

b) There must exist a subset $\frak h''\subseteq \frak h'$, $|\frak h''|=\rho$, $2\leq \rho\leq d$, and constants $\alpha_k\neq 0$, $k=1, 2, ..., \rho$, such that
  $$\sum_{k=1}^{\rho}\alpha_k\eta_k\equiv 0.\tag 2.26$$
Furthermore, (by lowering $\rho$ if necessary) one may assume that for all subsets $\sigma\subset \frak h''$, $|\sigma|=\rho-1$, 
$$\sum_{k=1}^{\rho-1}\beta_k\eta_k\equiv 0 \Longleftrightarrow  \text{ the constants } \beta_k=0, \text{ for all }k. \tag 2.27$$
This can be dealt with similarly to $\rho=d+1$ in case a1, ii),  as follows.

Using (2.26), (2.9) gives 
$$\Bbb Z^b\ni\sum_{k=1}^{\rho}\alpha_k\nu_k=0.$$
So 
$$\sum_{k=1}^{\rho}\alpha_k\nu_k\cdot\omega^{(0)}=0$$
and we have 
$$\sum_{k=1}^{\rho}\alpha_kJ_k=\sum_{k=1}^{\rho}\alpha_k|\eta_k|^2.$$
Solving for $\eta_{\rho}$ using (2.26), we have 
$$\eta_{\rho}=-\sum_{k=1}^{\rho-1}\frac{\alpha_k}{\alpha_\rho}\eta_k, \quad \alpha_k\neq 0, \,k=1, 2, ...,\rho.$$
So $$\sum_{k=1}^{\rho}\alpha_k\eta_k^2=\sum_{k\neq\rho}\alpha_k\eta_k^2+\frac{1}{\alpha_\rho}(\sum_{k\neq\rho}\alpha_k\eta_k)^2\not\equiv 0$$
when $\rho>2$, as the second sum yields cross terms. 

When $\rho=2$ and (2.26) holds, 
$$\sum_{k=1}^2\alpha_k\eta_k^2\equiv 0$$
if and only if $\alpha_1=-\alpha_2$, in this case $\eta_1\equiv \eta_2$, which is a contradiction. 

So when $$\sum_{k=1}^{\rho}\alpha_k\eta_k\equiv 0,\tag 2.28$$
$$\sum_{k=1}^{\rho}\alpha_kJ_k\not\equiv 0, \tag 2.29$$
for $2\leq \rho\leq d$.

Let $\pi$ be a projection of $\Bbb R^d$ to a $\Bbb R^{\rho-1}$ subspace and 
$$\zeta=\{\zeta_k\}_{k=1}^\rho=\{\pi \eta_k\}_{k=1}^{\rho},$$
where $\eta_k\in\frak h''$.  Define 
$$D_\zeta=\det[[ (2\zeta, \Cal J)]]_{\rho\times\rho},$$
where $\Cal J=\{J_k\}_{k=1}^\rho$.  Let $\Cal D_\zeta$ be the set 
defined by $D_\zeta=0$.  Using (2.27), there must exist 
a $\Bbb R^{\rho-1}$ subspace $Z$, $\zeta_k\in Z$, $k=1, 2, ..., \rho-1$, such that 
$$\det [[\{\zeta_k\}_{k=1}^{\rho-1}]]_{(\rho-1)\times(\rho-1)} \not\equiv 0.\tag 2.30$$
Determinant expansion using the $\rho$th row together with 
(2.28-2.30) then gives that 
$$D_\zeta\not\equiv 0.$$ 
So 
$$\Cal D'=\bigcap_\zeta \Cal D_\zeta$$
has co-dimension $1$,  
where the intersection is over all projections of $\Bbb R^d$ onto all possible $\Bbb R^{\rho-1}$ subspaces. 
On $(\Bbb R^d)^b\backslash \Cal D'$, (*) is satisfied. (Note that on $(\Bbb R^d)^b\backslash \Cal D'$, $\eta_k$, 
$k=1, 2, ..., \rho-1$, are linearly independent {\it as vectors} and {\it not} just as functions. This is because 
otherwise the determinant in (2.30) is $0$ for all $\zeta$, and therefore $D_\zeta=0$ for all $\zeta$, which 
is a contradiction.)  

By an abuse of notation, call all the sets defined by all the previous $D=0$, $\Cal D'$ as well. 
Take the union over all the possible $\Cal D'$ generated by all the possible subsets $\frak h$ 
(not necessarily connected) of $(d+1)$ elements of $\Cal W'$, that are included in no $L_k$, $k=1, 2, ..., b$, 
$L_k$ as defined in (Gi). Call the resulting set $\Cal D$. Then
$$\Cal D\supset\Upsilon_3$$
is algebraic of co-dimension $1$. The complement $\Cal D^c$ is Zariski open, and on it 
(*) is satisfied. This concludes the proof .
$\square$
\enddemo

\noindent{\it Remark.} In the proof of $\Upsilon_3$ being co-dimension $1$, the connected set property is not explicitly 
used, only that $\Upsilon_3$ is contained in the union of finite number of co-dimension $1$ sets determined by $\Cal W'$, which 
is an algebra to order $2(d+1)$. The latter restriction makes (Giii) useful mainly to bound the sizes of {\it connected} sets, cf. Lemma 2.2 and its
proof.
\smallskip
\demo{Proof of Lemma 2.3}
This is because $(\Bbb R^d)^b\backslash\Upsilon$ is Zariski open by the constructions in Lemmas 2.4 and 2.5
and therefore contains an infinite set of integers. \hfill $\square$ 
\enddemo

Below we briefly indicate the considerations that lead to (Gi-iii). For more details, see sect. 3. 

\subheading {2.3 Origins of the genericity conditions}

To implement the Newton scheme using (1.16), we need to bound $[F'_{\Cal S^c}]^{-1}$. From previous considerations,
it suffices to consider $[PA_0P]^{-1}$ with $P$ defined in (1.17), $A_0$ in (1.15). For $u^{(0)}$ satisfying (Gi-iii),  we show in sect. 3 that
$PA_0P=\oplus_\alpha\Cal A_\alpha$, where $\Cal A_\alpha$ are T\"oplitz matrices of sizes at most $(2b+2d)\times (2b+2d)$. 
This can be seen by using the notion of connected sets on $\Cal C$. 

Since we have made the identification: 
$$\Bbb Z^{b+d}\times\Bbb Z_2\sim\{\Bbb Z^{b+d},+\}\cup\{\Bbb Z^{b+d},-\},$$
we have
$$\Cal C=\{\Cal C^+,+\}\cup\{\Cal C^-,-\}.$$
For notational simplicity, we generally drop the $\pm$ signs and write $(n, j)\in\Cal C^+$ for $(n, j,+)\in\{\Cal C^+,+\}$
and $(n, j)\in\Cal C$ for either $(n, j,+)\in\{\Cal C^+,+\}$
or $(n, j,-)\in\{\Cal C^-,-\}$
etc.
\smallskip
\noindent{\bf Definition.} A subset $\frak s \subset \Cal C$ is {\it connected} if 
for all $a, b\in \frak s$, $\exists s_1, s_2, ..., s_k\in \frak s$, such that if we set $s_0=a$ and $s_{k+1}=b$, then 
$$\aligned  &s_{i+1}-s_i\in\Gamma\backslash \{(0, 0)\}, \text{  if } s_{i+1}, s_i\in\Cal C^+  \text{  or } s_{i+1}, s_i\in\Cal C^-;\\
&s_{i+1}-s_i\in\Gamma_+, \text{ if } s_{i+1}\in\Cal C^+, s_i\in\Cal C^-;\\
&s_{i+1}-s_i\in\Gamma_-, \text{ if } s_{i+1}\in\Cal C^-, s_i\in\Cal C^+;\endaligned\tag\dag
$$
for all $i=0, 1, 2, ..., k$, where we used the aforementioned convention of not making explicit the $\Bbb Z_2$ index.

\noindent{\it Remark.}  The above definition is just (2.21-2.22) restricted to $\Cal C^+\cup\Cal C^-$. 

For $$s_i=(n, j,+)\in\Cal C^+,$$
define 
$$-s_i:=(-n, -j,-)\in\Cal C^-.$$
Define $$-(-s_i)=s_i.$$ For $\frak s=\{s_i\}$, define the set $-\frak s$ to be 
$$-\frak s:=\{-s_i\}.$$ We note that if $\frak s$ is connected, then by reflection symmetry
$-\frak s$ is also connected. 

Every connected set $\frak s\subset\Cal C$ can be mapped to a connected set $\sigma\subset\Cal W$
by reversing the map defined in (2.18-2.20). The map is as follows:  

\noindent Assume $|\frak s|=R+1$ and $\frak s\cap \Cal C^+\neq\emptyset$. Let $s_1\in\frak s\cap\Cal C^+$. 
By abuse of notation, let $s_i$ also denote the $\Bbb Z^{b+d}$ component of $s_i$. Define 
$$\sigma'=\{s_k-s_i|k=1, 2, ..., R+1, s_k\in\frak s\}\subset\Bbb Z^{b+d}.$$ 
Then the set $\sigma=\sigma'\backslash\{(0, 0)\}$ is a connected set in $\Cal W_R$.
If $\frak s\cap \Cal C^+=\emptyset$, then $-\frak s\cap \Cal C^+\neq\emptyset$. The corresponding 
$-\sigma'$ gives a connected set $-\sigma$ in $\Cal W_R$.

We have so far defined connected sets on $\Bbb Z^{b+d}\times\Bbb Z_2$ as 
well as on the restriction $\Cal C$. The main reason for the latter is that a connected set $\frak s$ on $\Cal C$
provides additional equations that must be satisfied by elements of $\frak s$, see (2.31-2.33) below. This in turn limits the size $|\frak s|$
of $\frak s$ under the genericity conditions (Gi-iii). In sect. 3, we shall harvest the consequences of this restriction on size. 
Below we give an indication of the idea. 

Assume that $\frak s$ is a connected set on $\Cal C$ and that the sites 
to be mentioned are elements of $\frak s$. 
If $(n,j)\in\Cal C^+$ and $(n', j')\in\Cal C^+$ are connected, then 
$n'=n+\nu$ and $j'=j+\eta$, where $(\nu, \eta)\in\Cal W'$ and
$$
\cases
(n\cdot\omega^{(0)}+j^2)=0,\\
(n+\nu)\cdot\omega^{(0)}+(j+\eta)^2=0;
\endcases\tag 2.31
$$
and if $(n',j')\in\Cal C^-$, then $(\nu, \eta)\in\Cal W\backslash\Cal W'$ and
$$
\cases
(n\cdot\omega^{(0)}+j^2)=0,\\
-(n+\nu)\cdot\omega^{(0)}+(j+\eta)^2=0.
\endcases\tag 2.32
$$
Similar if  $(n_-,j_-)\in\Cal C^-$ and $(n'_-,  j'_-)\in\Cal C^-$ are connected then
$$
\cases
-n_-\cdot\omega^{(0)}+j^2_-=0,\\
-(n_-+\nu)\cdot\omega^{(0)}+(j_-+\eta)^2=0;
\endcases\tag 2.33 
$$
for some $(\nu, \eta)\in\Cal W'$. 

For each connected set $\frak s$ on $\Cal C$ described above, (2.31-2.33) define a system of polynomial equations. 
Up to reflection symmetry, we may assume that the system formed from (2.31) is at least as large as that from (2.33). 
Subtracting pairwise the second equations from the first equation in (2.31) and also pairwise (if there is at least one pair) the ones in (2.33)
lead to two systems of {\it linear} equations, each in $d$ variables, namely the $d$ components of $j$ or $j_-$. This gives rise to hyperplanes of types which appear in the genericity
conditions (Giii). We do not make use of the equations in (2.32), which upon addition yield quadratic equations in $j$ describing ellipsoids.

For $u^{(0)}$ satisfying (Gi-iii), 
we show by contradiction, in Lemma 3.2 of sect. 3, that the largest connected sets are of sizes at most $\text {max }(2b, 2d+2)\leq 2b+2d$. The 
exceptional connected sets of size 
$2b$ result from translation invariance and have spatial support the set $\{\pm j_k\}_{k=1}^b$. The other 
connected sets are of sizes at most $2d+2$. The translation invariant sets correspond to degeneracy and are
in fact the only reason for requiring the leading nonlinear $\Cal O(\delta^{2p+1})$ term in (1.1) to be independent of $x$. 
This is a sufficient but {\it not} necessary condition.
The $x$ dependence of the higher order terms does not matter as they are treated as perturbations.

\noindent{\it Remark.} It is possible that by making use of the quadratic equations in $j$, one could have a
better bound on the sizes of the connected sets away from the set $\{\pm j_k\}_{k=1}^b$. But a bound 
such as $2d+2$, {\it independent} of the number of frequencies $b$, suffices for the ensuing analysis
in sect. 4.  
\smallskip

\noindent{\it A note on the constants.}

Starting from sect. 3, there is a large number of positive constants which result from the estimates. The small constants are generally denoted
by $c$, $c'$ etc. and $\epsilon$, while the large ones by $C$, $C'$ etc. Unless indicated otherwise, they are not the same and may vary from 
statement to statement. In particular, the $\epsilon$'s that appear in sects. 3-5, are {\it not} the same as the $\epsilon$ in the Theorem in sect. 1.

\noindent {\it A note on the rescaling: $a\to\delta a$}

In sects. 3-5, we shall seek solutions with small amplitude $a=\{a_k\}_{k=1}^b$, and therefore rescale:
 $a_k\to\delta a_k$ ($0<\delta\ll 1$) for $k=1, ..., b$. (In the statement of the Theorem, we return to the ``original $a$".)
\bigskip 

 \head {\bf 3. The first step in the Newton scheme -- extraction of parameters} \endhead
 Starting in this section,  we write $F'$ for  the $F'_{\Cal S^c}$ in (1.16), as it will only appear in the context of  the $P$-equations; similarly 
 we write $F$ for the $F_{\Cal S^c}$ in (1.16). (Recall that $F|_{\Cal S}=0$ always.) In Lemma 3.1, we prove $F'$ is invertible with exponential off-diagonal decay.
The proof  rests on Lemma 3.2, which shows that under the genericity conditions the operator 
 $PA_0P$ as defined by (1.17, 1.15) is a block diagonal matrix. 
 
As a consequence, ${F'}^{-1}$ is controlled by a finite family of polynomials in $a=\{a_k\}_{k=1}^b$ -- the 
 amplitude -- the Fourier coefficients of the unperturbed solution $u^{(0)}$. These polynomials are the determinants of the 
 block matrices in the direct sum decomposition of $PA_0P$.  Varying $a$ thus leads to invertibility of $F'$ on {\it open} sets; moreover, pairing with the resolvent expansion yields exponential off-diagonal decay of
 ${F'}^{-1}$. 
 
 In Proposition 3.3, the $P$-equations produce the first corrections: 
 $$\Delta u^{(1)}=u^{(1)}-u^{(0)},\text{ and } \Delta v^{(1)}=v^{(1)}-v^{(0)};$$ 
 solving the $Q$-equations gives the modulated frequencies $\omega^{(1)}$, and that 
 the amplitude-frequency map 
 $a\to\omega^{(1)}(a)$ is a diffeomorphism.
 
 \subheading {3.1 The invertibility of $F'$} 
 
After rescaling $a\to\delta a$, we solve instead:
$$
\cases
\text{diag }(n\cdot\omega+j^2)u+\delta^{2p} (u*v)^{*p}* u+\sum_{m=1}^\infty\delta^{2p+2m}\alpha_m*(u*v)^{*(p+m)}*u=0,\\
\text{diag }(-n\cdot\omega+j^2)v+\delta^{2p} (u*v)^{*p}* v+\sum_{m=1}^\infty\delta^{2p+2m}\alpha_m*(u*v)^{*(p+m)}*v=0,
\endcases\tag 3.1
$$
with $u|_{\text{supp } u^{(0)}}=u^{(0)}=a\in (0,1]^b=\Cal B(0,1)$, and similarly for $v$. 

\proclaim{Lemma 3.1} 
Assume that $u^{(0)}=\sum_{k=1}^b a_k e^{-ij_k^{2}t}e^{ij_k\cdot x}$ a solution to the linear equation with $b$ frequencies 
satisfies genericity conditions (Gi-iii)
and $a=\{a_k\}\in (0, 1]^b=\Cal B(0,1)=\Cal B\subset\Bbb R^b\backslash\{0\}$. There exist $C$, $c>0$, such that
for all $\epsilon\in (0,1)$, there exists $\delta_0>0$ and for all $\delta\in (0, \delta_0)$, a set $\Cal B' _{\epsilon, \delta}\subset\Cal B$ with 
$$\text{meas }\Cal B'_{\epsilon, \delta} <C\delta^{c\epsilon}.$$
If $a\in\Cal B\backslash\Cal B' _{\epsilon,\delta}$, then 
$$\Vert [F'(u^{(0)}, v^{(0)})]^{-1}\Vert\leq \Cal O(\delta^{-2p-\epsilon}).\tag 3.2$$
Let $\pi$ be the projection of $\Bbb Z^{b+d}\times\Bbb Z_2$ onto $\Bbb Z^{b+d}$. 
There exists $\beta\in (0,1)$ such that
$$|[F'(u^{(0)}, v^{(0)})]^{-1}(x,y)|\leq \delta^{\beta|x-y|}=e^{-\beta|\log\delta||x-y|}, \tag 3.3$$
for all $|x-y|>1/\beta^2$, where $x$, $y\in\Bbb Z^{b+d}\times\Bbb Z_2$
and 
$$|x-y|:=|\pi x-\pi y|.\tag 3.4$$
\endproclaim 
\smallskip
As before, let $P_\pm$ be the projection on $\Bbb Z^{b+d}$ onto $\Cal C^\pm$ defined in (1.5, 1.6), and
$$P=\pmatrix P_+&0\\0&P_-\endpmatrix\tag 3.5$$
on $\Bbb Z^{b+d}\times \Bbb Z_2$. Let $A_0$ be as defined in (1.15) and 
$$A':=PA_0P\tag 3.6$$
be the restricted operator on $\Cal C$. The proof of Lemma 3.1 rests on the following geometric 
operator decomposition: 
\proclaim{Lemma 3.2}
Assume that $u^{(0)}=\sum_{k=1}^b a_k e^{-ij_k^{2}t}e^{ij_k\cdot x}$ a solution to the linear equation with $b$ frequencies 
satisfies the genericity conditions (Gi-iii). Then $A'$ can be written as 
$$A'=A'(a)=\oplus A'_{\alpha} (a),\tag 3.7$$
where $\alpha$ are connected sets on $\Cal C$ defined as in (\dag), satisfying 
$|\alpha|\leq 2b+2d$, $A'_{\alpha}$ are $A'$ restricted to $\alpha$ and therefore matrices of sizes at most $(2b+2d)\times (2b+2d)$.

Let $\pi$ be the projection of $\Bbb Z^{b+d}$ onto $\Bbb Z^d$.  The connected sets $\alpha$ have the further characterizations:
$$\align&\text{if } \pi\alpha\subseteq\{j_k\}_{k=1}^b\cup \{-j_k\}_{k=1}^b,\\
&\text{then } \pi\alpha=\{j_k\}_{k=1}^b\cup \{-j_k\}_{k=1}^b, \text{so } |\alpha|=2b,\tag 3.8 \\
&\text{if }\pi\alpha\not\subseteq\{j_k\}_{k=1}^b\cup \{-j_k\}_{k=1}^b, \text{ then } |\alpha|\leq 2d+2, \tag 3.9\endalign$$
where for notational simplicity $\{\pm j_k\}_{k=1}^b$ denote $\{\pm j_k, \pm\}_{k=1}^b$.
\endproclaim
\demo{Proof}
Assume that there is a connected set $\alpha$ on $\Cal C$ of length $|\alpha|=2d+3$. $\alpha$ can be written as $\alpha=\alpha^+\cup \alpha^-$ with $\alpha^+\subset\Cal C^+$
and $\alpha^-\subset\Cal C^-$. 

Without loss of generality, we may assume $|\alpha^+|\geq |\alpha^-|$. So $|\alpha^+|\geq d+2$, and
there must exist ${\tilde \alpha}^+\subseteq \alpha^+$, $|{\tilde \alpha}^+|=d+2$, such that ${\tilde \alpha}^+\cup \alpha^-$ is a connected set.  (The other case works the same way using reflection symmetry.)
Assume $(n,j)\in {\tilde \alpha}^+$, the set ${\tilde \alpha}^+$ then gives a system of $|{\tilde \alpha}^+|$ quadratic (in $j$) polynomials of the form:
$$
\cases
(n\cdot\omega^{(0)}+j^2)=0,\qquad\qquad\qquad\qquad\qquad\qquad\qquad\qquad\qquad\qquad\qquad\quad\,\,\,(3.10)\\
(n+\nu)\cdot\omega^{(0)}+(j+\eta)^2=0;\qquad\qquad\qquad\qquad\qquad\qquad\qquad\qquad\qquad\quad\,(3.11)
\endcases
$$
where $(\nu, \eta)\in\Cal W'$.

Subtracting pairwise the equations in (3.11) from (3.10), when (Gii) is satisfied, we obtain a system of $|{\tilde \alpha}^+|-1=d+1$ linear equations 
of the form: 
$$\align&2j\cdot\eta+\eta^2+\nu\cdot\omega^{(0)}\\
=\,&2j\cdot\eta+J=0,\endalign$$
where $J$ as defined in (2.24). (When (Gii) is violated, there could be less than $d+1$ equations, if $\nu\cdot\omega^{(0)}=0$, when $\eta=0$.)
There are $2$ possibilities.

\noindent {a)} The conditions on $\frak{h}$ in (Giii) are not violated,  then (Giii) gives that
$|{\tilde \alpha}^+|\leq d+1$, which is a contradiction to $|{\tilde \alpha}^+|=d+2$. Therefore $|\alpha^+|\leq d+1$ and $|\alpha|\leq 2|\alpha^+|\leq 2(d+1)$.

\noindent {b)} The conditions on $\frak{h}$ in (Giii) are violated. In this case, $b$ must satisfy $b>d+1$. This is because the exceptional sets $L_k$ in (Giii) satisfy
$|L_k|=b-1$. So if $b\leq d+1$, then $|L_k|\leq d$. Hence $\frak h\not\subseteq L_k$, as $|\frak h|=d+1$.  Therefore there must exist $k$ such that 
$$(2\frak h, \Cal J)\subseteq \{(2(j_i-j_k), -2(j_i-j_k)\cdot j_k), i=1, 2, ..., b, i\neq k\}$$
and 
$$D=\text{ det}[[(2\frak{h}, \Cal J)]]_{(d+1)\times(d+1)}\equiv 0.$$
The solutions to the equations in (3.10, 3.11) 
have $j$ coordinates in $\{ j_k\}_{k=1}^b$ and are the only solutions using (Gi).  The set $\alpha$ is a subset of a 
{\it maximally} connected set $\bar \alpha$, i.e., if $(n', j')$ is connected 
to $(n, j)\in\bar \alpha$, then  $(n', j')\in\bar \alpha$, satisfying 
$\pi\bar \alpha=\{j_k\}_{k=1}^b\cup \{-j_k\}_{k=1}^b$ and $|\bar \alpha|=2b$ using also (Gii).

This can be checked directly.  It follows directly from the definition in (2.21-2.22) that
$$E_\mu=\{(-e_k+\mu, j_k)\}_{k=1}^b\cup \{(e_k+\mu, -j_k)\}_{k=1}^b , \tag $\diamondsuit$  $$
$ \mu\in \Bbb Z^d, \mu\cdot\omega^{(0)}=0$
are connected sets. The projection 
$$\pi E_\mu= \{j_k\}_{k=1}^b\cup \{-j_k\}_{k=1}^b.$$ 
 On $E_\mu$, the condition in (Giii) is violated; moreover using (Gi), these are 
the only sets that violate the condition. 

We first prove that $E_\mu$ and $E_{\mu'}$ are {\it not} connected if $\mu\neq \mu'$, $\mu, \mu'\in\Bbb Z^d$, i.e., 
$\not\exists (n, j_i)\in E_\mu$, $(n', j_k)\in E_{\mu'}$, $(n'', -j_m)\in E_{\mu'}$, such that 
either 
$$(n'-n, j_k-j_i)\in\Gamma\backslash\{(0, 0)\}, \tag  3.12$$
or  $$(n''-n, -j_m-j_i)\in\Gamma_-. \tag 3.13$$
(Using reflection symmetry, it suffices to consider the above two cases. ) 

If (3.12) holds, then $$(n'-n, j_k-j_i)=(e_i-e_k+\mu'-\mu, j_k-j_i)\in \Gamma\backslash\{(0, 0)\}.$$
Since $$(-e_i+e_k, - j_k+j_i)\in \Gamma,$$ this shows that 
$$(\mu'-\mu, 0)\in 2\Gamma$$
for $\mu-\mu'\neq 0$, which contradicts (Gii), since $2\Gamma\backslash\{(0, 0)\}\subset \Cal W'$.

If (3.13) holds, 
then $$(n''-n, -j_m-j_i)=(e_m+e_i+\mu'-\mu, -j_m-j_i)\in \Gamma_-.$$
Since $$(-e_m-e_i,  j_m+j_i)\in \Gamma_+,$$ this shows that 
$$(\mu'-\mu, 0)\in \Gamma_++\Gamma_-=2\Gamma$$
for $\mu-\mu'\neq 0$, again contradicting (Gii). 

We now prove that $E_\mu$ are {\it maximal}. Assume that $(n, j)\in \Cal C$ is such that 
$$j\not\in\{j_k\}_{k=1}^b\cup \{-j_k\}_{k=1}^b.$$
Without loss of generality, one may assume 
$(n, j)\in\Cal C^+$, then for {\it no} $k$, 
$$j-j_k\in L_k.$$
So for {\it no} $\mu$
is $$S_\mu=\{\{(n, j)\}\cup E_\mu\}$$
connected. This is because if $S_\mu$ is connected, then 
$|S_\mu|\leq 2(d+1)$ from (Giii). But 
$|S_\mu|=|E_\mu|+1=2b+1>2(d+1)$ for $b>d+1$, since $E_\mu$ violates the condition in (Giii). 
So $E_\mu$ are maximal for all $\mu$.
This proves that $$\alpha\subseteq E_\mu:=\bar \alpha$$
for some $\mu$. 
\smallskip

Take $\alpha$ to be {\it maximal} connected sets on $\Cal C$. The conclusion in b) 
proves the property (3.8); while a) together with (Gi) prove (3.9). The direct sum decomposition in (3.7)
follows from this geometric decomposition of $\Cal C$ by $\alpha$. 
\hfill$\square$
\enddemo

\demo{Proof of Lemma 3.1}  (i) {\it The norm estimates}

Let $P^c=\Bbb I-P$, where $P$ as defined in (3.5) . The linearized operator  $F'(u^{(0)}, v^{(0)})$
is $F'=D+A$ with $D$ and $A=\delta^{2p}A_0+\Cal O(\delta^{2p+2})$ as in (1.14, 1.15). 

The Schur complement reduction \cite{S1, 2} implies that, $\lambda\in\sigma(F')\cap [-1/2,1/2]$ if and only if
$0\in\sigma(\Cal H)$, where 
$$\Cal H=PF'P-\lambda+PF'P^c(P^cF'P^c-\lambda)^{-1}P^cF'P$$
is the effective operator acting on the bi-characteristics $\Cal C$. 
Since $\Vert P^cF'P^c\Vert>1-\Cal O(\delta^{2p})>1/2$ and $\Vert PF'P^c\Vert=\Cal O(\delta^{2p})$, 
the last term is of order $\delta^{4p}$, uniformly for $\lambda\in[-1/2, 1/2]$. 

So $$\Cal H=PF_0'P-\lambda+\Cal O(\delta^{2p+2})\tag 3.14$$ in $L^2$ uniformly for $\lambda\in[-1/2, 1/2]$,
where $$F_0=\delta^{2p}\pmatrix (u*v)^{*p}*u\\(u*v)^{*p}*v\endpmatrix. \tag 3.15$$
To obtain (3.2),
it suffices to prove $\Vert [PF_0'P]^{-1}\Vert\leq \Cal O(\delta^{-2p-\epsilon})$. It is important to note that since 
$$PF_0'P=\delta^{2p} PA_0P=\delta^{2p}A',$$ where $A_0$, $A'$ as in (1.15, 3.6),
$A'=A'(a)$ depends on $a$ but is {\it independent} of $\delta$. 

From Lemma 3.2 
$$A'(a)=\oplus A'_{\alpha} (a),$$
where $\alpha$ are connected sets on $\Cal C$ and $A'_{\alpha}$ are matrices of sizes at most $(2b+2d)\times (2b+2d)$.
Since $A_0(a)$ as defined in (1.15)  is a convolution matrix on $\Bbb Z^{b+d}\times\Bbb Z_2$ and $A_0(x, y)=A_0(x-y, 0)\neq 0$ for at most  $2b^{2p}$ of $(x-y)$, where $x$, $y\in\Bbb Z^{b+d}\times\Bbb Z_2$, there are at most $2^{2b^{2p}}=4^{b^{2p}}=K$ types of $A'_{\alpha}$ using that $A'$ is the restriction of $A_0$ to $\Cal C$ as defined in (3.6).
We rename $A'_{\alpha}$ as $\Cal A_k$, $1\leq k\leq K$. 

For each $k$, $\det \Cal A_k(a)=P_k(a)$ is a polynomial in $a$ of degree at most 
$2p(2b+2d)$. For all $k$, $P_k$ is a non-constant function on $\Cal B=(0,1]^b$, which can be seen as follows.
Set $a=(a_1, 0,..,0)\in\Bbb R^b$. The convolution matrices $u*v$ and $u*u$ are then with matrix elements:
$$(u*v)[(n, j), (n, j)]=|a_1|^2=a_1^2, (n, j)\in\Bbb Z^{b+d},  \text{ and } 0 \text{ otherwise};$$
and   
$$(u*u)[(n, j), (n', j')]=a_1^2 \text{ if } (n, j)-(n', j')=(-2e_1, 2j_1)  \text{ and } 0 \text{ otherwise},$$
where $e_1:=(1, 0,..., 0)\in\Bbb Z^b$ as before and $j_1\in\Bbb Z^d$.

To see the structure of $A_0$ we rewrite $A_0$ as before and have 
$$A_0=(A_0)_{11}\otimes \pmatrix 1&0\\0&0\endpmatrix+ (A_0)_{12}\otimes \pmatrix 0&1\\0&0\endpmatrix+ (A_0)_{21}\otimes \pmatrix 0&0\\1&0\endpmatrix+(A_0)_{22}\otimes \pmatrix 0&0\\0&1\endpmatrix,$$
where $(A_0)_{11}=(A_0)_{22}$ are diagonal matrices operating on $\ell^2(\Bbb Z^{b+d})$ with all diagonal 
elements equal to $$(p+1)|a_1|^{2p}=(p+1)a_1^{2p},$$ $(A_0)_{12}$ is the convolution matrix on $\ell^2(\Bbb Z^{b+d})$
with matrix elements
$$(A_0)_{12}[(n, j), (n+2e_1, j-2j_1)]=p|a_1|^{2(p-1)}a_1^2=pa_1^{2p}\text{ and } 0\text{ otherwise},$$
and $(A_0)_{21}=(A_0)_{12}^{\text t}$ is the transpose.

Taking the decomposition in (3.7) into account, $\Cal A_k$ is then either a diagonal matrix with diagonal elements 
$(p+1)a_1^{2p}$ or a matrix with diagonal elements $(p+1)a_1^{2p}$ and $2$ non-zero off diagonal elements
both equal to $pa_1^{2p}$. The determinant $P_k$ is therefore not a constant. So there exist $C$, $c>0$, such that for all 
$0<\epsilon<1$, there exists $\delta_0>0$, such that for all $\delta\in (0, \delta_0)$, 
$$\text{ meas }\{a\in \Cal B||P_k|<\delta^\epsilon, \text{ all }k\leq K\}\leq C\delta^{c\epsilon},\tag 3.16$$
which is a simple consequence of polynomial $P_k$. However the corresponding estimates  
hold in the general analytic category, cf. e.g., \cite{Lemma 11.4, GS} and references therein.

Since $\Vert \Cal A_k\Vert\leq \Cal O(1)$, if $|\det\Cal A_k|>\delta^\epsilon$, then 
$\Vert [\Cal A_k]^{-1}\Vert\leq\Cal O(\delta^{-\epsilon})$. (The exponent is $-1$ because of self-adjointness.)
In view of (3.14), this proves (3.2).

\noindent (ii) {\it The point-wise estimates}

To prove (3.3), let
$$\align \tilde F&=\oplus_\alpha \delta^{2p}A'_{\alpha}\oplus[\text{ diag }(n\cdot\omega^{(0)}+j^2+\delta^{2p}(A_0)_{11}(0, 0))|_{(n,j)\notin\Cal C^+}\otimes\pmatrix 1&0\\0&0\endpmatrix]\\
&\qquad\oplus[\text{ diag }(-n'\cdot\omega^{(0)}+{j'}^2+\delta^{2p}(A_0)_{22}(0, 0))_{(n',j')\notin\Cal C^-}\otimes\pmatrix 0&0\\0&1\endpmatrix],\\
:&=\oplus_\alpha \delta^{2p}A'_{\alpha}\oplus \Cal D_{\Bbb Z^{b+d}\times \Bbb Z_2\backslash\Cal C}\tag 3.17\endalign$$
where the first direct sum is exactly as in (3.7), with $\alpha$ connected subsets of $\Cal C$ and $(A_0)_{ii}(0, 0)$ denotes the diagonal element of the convolution matrix $(A_0)_{ii}$, $i=1, 2$. 
Using (3.17), the resolvent expansion gives:
$$[F']^{-1}=[\tilde F]^{-1}-[\tilde F]^{-1}\tilde\Gamma[\tilde F]^{-1}+[\tilde F]^{-1}\tilde\Gamma[\tilde F]^{-1}\tilde\Gamma[F']^{-1},\tag 3.18$$
where 
$$\tilde\Gamma=F'-\tilde F:=\Gamma_1+\Gamma_2,\tag 3.19$$
satisfying $$\Vert\tilde\Gamma\Vert=\Cal O(\delta^{2p});\tag 3.20$$
$$\Gamma_1=(D+\delta^{2p}A_0)-\tilde F$$
with $D$, $A_0$ as in (1.14, 1.15),
and 
$\Gamma_2$ is the convolution operator 
$$\Gamma_2=A-\delta^{2p}A_0$$
with $A$ as defined in (1.15). So 
$$\Vert\Gamma_2\Vert=\Cal O(\delta^{2p+2}).\tag 3.21$$

Since  $$\Vert \tilde F^{-1}\Vert=\Cal O(\delta^{-2p-\epsilon}),\tag 3.22$$
from (3.17), using (3.20), we have $$\Vert \tilde F^{-1}\tilde\Gamma\Vert=\Cal O(\delta^{-\epsilon}).\tag 3.23$$
To estimate $$\Vert \tilde F^{-1}\tilde\Gamma\tilde F^{-1}\tilde \Gamma\Vert,$$
we write
$$\aligned\tilde F^{-1}\tilde\Gamma\tilde F^{-1}\tilde\Gamma&=\tilde F^{-1}\Gamma_1\tilde F^{-1}\Gamma_1+\tilde F^{-1}\Gamma_1\tilde F^{-1}\Gamma_2+\tilde F^{-1}\Gamma_2\tilde F^{-1}\Gamma_1+\tilde F^{-1}\Gamma_2\tilde F^{-1}\Gamma_2\\
:&=\tilde F^{-1}\Gamma_1\tilde F^{-1}\Gamma_1+\Cal O(\delta^{2-2\epsilon}),\endaligned\tag 3.24$$
where we used (3.19-3.22) to estimate the last three terms and $\Cal O(\delta^{2-2\epsilon})$ is in operator norm. 

To estimate the first term, we use (3.17) and write $\Cal D$ for  $\Cal D_{\Bbb Z^{b+d}\times \Bbb Z_2\backslash\Cal C}$. We have
$$\aligned \tilde F^{-1}\Gamma_1\tilde F^{-1}= &[\oplus_\alpha \delta^{-2p}{A'}^{-1}_{\alpha}]\Gamma_1 [\oplus_{\alpha'} \delta^{-2p}{A'}^{-1}_{\alpha'}]+[\oplus_\alpha  \delta^{-2p}{A'}^{-1}_{\alpha}]\Gamma_1 [\Cal D^{-1}]\\
&+[\Cal D^{-1}]\Gamma_1 [\oplus_{\alpha'}  \delta^{-2p}{A'}^{-1}_{\alpha'}]+[\Cal D^{-1}]\Gamma_1 [\Cal D^{-1}].\endaligned\tag 3.25$$
The first term in (3.25) is identically $0$. This is because from the definition of $\Gamma_1$, 
$$\Gamma_1(x, y)=\Gamma_1(y, x)=(A_0-\Cal A I)(x,y),$$
where $\Cal A=A_0(x, x)$, $I$ is the identity matrix, if $x\in\Cal C$, $y\in\Bbb Z^{b+d}\backslash \Cal C$;
and $$\Gamma_1(x, y)=0$$ otherwise.

Using $\Vert\Cal D^{-1}\Vert =\Cal O(1)$ in (3.25), we obtain 
$$\Vert \tilde F^{-1}\Gamma_1\tilde F^{-1}\Vert=\Cal O(\delta^{-\epsilon}).$$
Using the above in (3.24) and since $\Vert \Gamma_1\Vert=\Cal O(\delta^{2p})$, we obtain 
$$\Vert \tilde F^{-1}\tilde\Gamma\tilde F^{-1}\tilde\Gamma\Vert \leq\Vert \tilde F^{-1}\Gamma_1\tilde F^{-1}\Gamma_1\Vert+\Cal O(\delta^{2-2\epsilon})=\Cal O(\delta^{2-2\epsilon}).\tag 3.26$$
So combining (3.26) and (3.23), we have 
$$\Vert [\tilde F^{-1}\tilde\Gamma]^{2m}]\Vert\leq\Cal O( \delta^{2m(1-\epsilon)}),\, m=1, 2, ...\tag 3.27$$
and  
$$\Vert [\tilde F^{-1}\tilde\Gamma]^{2m-1}]\Vert\leq \Cal O(\delta^{2(m-1)(1-\epsilon)-\epsilon}),\, m=1, 2, ...\tag 3.28$$

Iterating the resolvent expansion in (3.18) $r$ times yields the $(r+1)$ term series 
$$[F']^{-1}=[\tilde F]^{-1}-[\tilde F]^{-1}\tilde\Gamma[\tilde F]^{-1}+...+(-1)^r[{\tilde F}^{-1}\tilde\Gamma]^r[F']^{-1}.$$
We note that the blocks in $\tilde F$ (and hence $[\tilde F]^{-1}$) are of sizes at most $(2b+2d)$ and that the symbols
of $\tilde \Gamma$ are trigonometric polynomials. Matrix multiplication then infers that for 
some $\beta>0$ depending only on $\text{supp } u^{(0)}$, $p$, $b$, $d$ and $H$ in (1.1), and any given $x$, $y\in\Bbb Z^{b+d}\times\Bbb Z_2$ 
satisfying $|x-y|>1/\beta^2$, there exists $r> 1$ such that the first $r$ terms in the series are identically $0$.
Using the bounds in (3.27, 3.28, 3.2) to estimate the last, the $(r+1)$th term 
produces (3.3).
\hfill $\square$
\enddemo

\noindent {\it Remark.} For simplicity of exposition, the last part of the proof uses the finite range nature of $\tilde\Gamma$. This is in fact not needed -- the 
result holds for $\tilde \Gamma$ with analytic symbols using weighted estimates, see e.g., \cite{W1}.

\subheading {3.2 The first iteration}

Using Lemma 3.1 to solve the $P$, and then the $Q$-equations, we obtain the following 
result after the first iteration. As earlier, let 
$$\Delta u^{(1)}=u^{(1)}-u^{(0)},\, \Delta v^{(1)}=v^{(1)}-v^{(0)}, \Delta \omega^{(1)}=\omega^{(1)}-\omega^{(0)}.$$
\proclaim{Proposition 3.3}  Assume that $u^{(0)}=\sum_{k=1}^b a_k e^{-ij_k^{2}t}e^{ij_k\cdot x}$ a solution to the linear equation with $b$ frequencies is generic and $a=\{a_k\}\in (0, 1]^b=\Cal B(0, 1)=\Cal B\subset\Bbb R^b\backslash\{0\}$. 
Let $\epsilon, \epsilon'\in (0, 1)$. There exists 
$\delta_0>0$, such that for all $\delta\in (0, \delta_0)$, there is a set $\Cal B_{\epsilon, \delta}$, $\Cal B\supset \Cal B_{\epsilon, \delta}\supset \Cal B'_{\epsilon, \delta}$ (the set in Lemma 3.1) with 
$$\text{meas }\Cal B_{\epsilon, \delta} <\epsilon'/2.$$ 
Let $\rho$ be a weight on $\Bbb Z^{b+d}$ satisfying 
$$\aligned\rho(x)=&e^{\beta|\log\delta| |x|},\, 0<\beta<1\, \text{ for } |x|>1/\beta^2,\\ 
=&1,\qquad\qquad\quad\,\qquad \quad\,\text{  for } |x|\leq 1/\beta^2.\endaligned$$
Define the weighted $\ell^2$ norm:
$$\Vert\cdot \Vert_{\ell^2(\rho)}=\Vert\rho\cdot\Vert_{\ell^2}.$$

There exists $\beta\in(0, 1)$, determined only by  
$\text{supp } u^{(0)}$, $p$, $b$, $d$ and $H$ in (1.1), such that if $a\in\Cal B\backslash\Cal B_{\epsilon, \delta}$,
then 
$$\align&\Vert \Delta u^{(1)}\Vert _{\ell^2(\rho)}=\Vert \Delta v^{(1)}\Vert_{\ell^2(\rho)}= \Cal O(\delta^{3-\epsilon}),\tag 3.29\\
&\Vert F(u^{(1)}, v^{(1)})\Vert_{\ell^2(\rho)\times\ell^2(\rho)}=\Cal O(\delta^{2p+5-2\epsilon}),\tag 3.30\\
&\Vert \Delta \omega^{(1)}\Vert\asymp \delta^{2p},\tag 3.31\\
&\Vert \frac{\partial \omega^{(1)}}{\partial a}\Vert\asymp \delta^{2p},\tag 3.32\\
&\Vert(\frac{\partial \omega^{(1)}}{\partial a})^{-1}\Vert\lesssim\Cal O_{\epsilon'}( \delta^{-2p}),\tag 3.33\\
&\big| \det(\frac{\partial \omega^{(1)}}{\partial a})\big|\gtrsim \Cal O_{\epsilon'}(\delta^{2pb}).\tag 3.34\endalign$$
(Note that for the first iteration, $\omega^{(1)}$ is defined on $\Cal B(0, 1)$, cf. (1.18), so 
$\frac{\partial \omega^{(1)}}{\partial a}$ is meant in the classical sense.)

Moreover $\omega^{(1)}$ is Diophantine
$$\Vert n\cdot \omega^{(1)}\Vert_{\Bbb T}\geq\frac{\kappa\delta^{2p}}{|n|^\gamma},\quad n\in\Bbb Z^b\backslash\{0\},\,\kappa>0, \gamma>2b+1,\tag 3.35$$
where $\Vert\,\Vert_{\Bbb T}$ denotes the distance to integers in $\Bbb R$, $\kappa$ and $\gamma$ are independent of $\delta$.
\endproclaim

\demo{Proof}
 {(i) {\it The $P$-equations} 
 
From the Newton scheme
$$\aligned \Delta \pmatrix u^{(1)}\\ v^{(1)}\endpmatrix&=[F'(u^{(0)}, v^{(0)})]^{-1}F(u^{(0)}, v^{(0)})\\
&=[{D'}^{-1}+{F'}^{-1}(A-\delta^{2p}\text{diag }A_0){D'}^{-1}]F,\endaligned$$
where $\text{diag }A_0$ is the diagonal part of $A_0$, $D'=D+\delta^{2p}\text{diag }A_0$, $D$ and $A_0$ as in (1.14, 1.15).
The matrix elements of $D'$ satisfy
$$|D'(n, j; n, j)|\geq\Cal O(\delta^{2p}), $$
for all $(n, j)$, since $D(n, j; n, j)\in\Bbb Z$, 
so $$\Vert {D'}^{-1}\Vert \leq \Cal O(\delta^{-2p}).$$ 

Let $F_0(u^{(0)}, v^{(0)})$ be as in (3.15), $(\nu, \eta)\in\Gamma\backslash\{(0, 0)\}$ and 
$(-e_m, j_m)\in \text{supp }u^{(0)}$. If $(\nu, \eta)=(-e_{k'}+e_m, j_{k'}-j_m)$ for some $k'=1,...,b$, 
$$(-e_m, j_m)+(\nu, \eta)=(-e_{k'}, j_{k'})\in\Cal S.$$
Otherwise (Giv) gives  
$$(-e_m, j_m)+(\nu, \eta)\notin\Cal C.$$
Therefore since 
$$\text{supp }F_0 (u^{(0)}, v^{(0)})=\text{supp } (u^{(0)}*v^{(0)})^{*p}*u^{(0)}\bigcup \text{supp } (u^{(0)}*v^{(0)})^{*p}*v^{(0)},$$
and $$\text{supp } (u^{(0)}*v^{(0)})^{*p}=\Gamma$$
from (2.8), (Giv) yields 
$$\text{supp }F_0(u^{(0)}, v^{(0)})\cap\{\Cal C\backslash\Cal S\}=\emptyset.$$
So 
$$\aligned \Vert {D'}^{-1}F\Vert_{\ell^2}&=\Vert {D'}^{-1}F_0+{D'}^{-1}(F-F_0)\Vert _{\ell^2}\\
&\leq \Cal O(\delta^{2p+1})+\Cal O(\delta^3)\\
&=\Cal O(\delta^3),\endaligned$$
and 
$$\Vert {F'}^{-1}(A-\delta^{2p}\text{diag }A_0){D'}^{-1}F\Vert_{\ell^2}=\Cal O(\delta^{3-\epsilon}),$$
where we also used (3.2).
So $$\Vert  \Delta \pmatrix u^{(1)}\\ v^{(1)}\endpmatrix\Vert_{\ell^2}=\Cal O(\delta^{3-\epsilon})$$
and $$\Vert F(u^{(1)}, v^{(1)})\Vert_{\ell^2}\leq \Cal O(\Vert F''\Vert)\Vert  \Delta \pmatrix u^{(1)}\\ v^{(1)}\endpmatrix\Vert^2_{\ell^2}=\Cal O(\delta^{2p+5-2\epsilon}).$$
Using (3.3), the above two estimates hold in weighted space as well and we obtain (3.29, 3.30). 

\noindent (ii) {\it The $Q$-equations} 

From the $Q$ equations, the frequency modulation is as before  
$$\aligned\Delta \omega_k^{(1)}&=\frac{1}{a_k}F(u^{(0)}, v^{(0)})(-e_k,j_k)=\delta^{2p}\frac{{(u^{(0)}*v^{(0)}})^{*p}*u^{(0)}}{a_k}(-e_k, j_k)+\Cal O(\delta^{2p+2}),\\
:&=\delta^{2p}\Omega_k+\Cal O(\delta^{2p+2}),\, k=1, 2, ..., b. \endaligned\tag 3.36$$
Since $a_i>0$, $i=1, 2, ..., b$, and $\Delta \omega_k^{(1)}$ are finite sums of 
polynomials in $a$, (3.31, 3.32) are immediate. (Note that the denominator must contain $a_k$ because of the restriction 
to $(-e_k, j_k)$, so the right hand side of (3.36) represents polynomials.)


To prove (3.33, 3.34), set the derivative matrix 
$$\frac{\partial \Omega}{\partial a}:=\big [\big[\frac{\partial\Omega_k}{\partial a_i}\big]\big], k, i=1, ...,b.$$ 
To see its structure, it suffices to analyze the 
polynomials given by $\Omega_k$. 

Let $M$ be the convolution matrix  ${(u^{(0)}*v^{(0)}})^{*p}*$. (Note that $M=(A_0)_{11}/{p+1}$ from (1.15).) 
Using (3.36), we have 
$$\Omega_k=M_{kk}+\sum_{i\neq k}\frac{M_{ki} a_i}{a_k}, \, k=1, ..., b.$$ 
From the structure of $M$, 
$$M_{kk}(a_1, a_2, ..., a_b)=P(a_1, a_2, ..., a_b)$$
and $$M_{ki}=P'(a_1, a_2, ..., a_b)a_k\bar a_i,$$
where $P$ and $P'$ are homogeneous polynomials in $a$ with positive integer coefficients
and are invariant under any permutations of the arguments, $P$ is of degree $2p$, $P'$, $2(p-1)$.
So $\Omega_k$ is a homogeneous polynomial in $\{a_i\}_{i=1}^b$ of degree $2p$ and can be written as 
$$\Omega_k(a_1, a_2, ..., a_b)=P(a_1, a_2, ..., a_b)+P'(a_1, a_2, ..., a_b)P_k(\{a_i\})_{i\neq k},$$
with $P_k=\sum_{i\neq k} a_i^2$ and $P$, $P'$ as above.

Set $a=(1,1, .., 1)$, we therefore have
$$\frac{\partial \Omega_k}{\partial a_i}(1, 1, ..., 1)> \frac{\partial \Omega_k}{\partial a_k}(1, 1, ..., 1)$$
for all $i\neq k$. Let $q$ be the diagonal elements and $Q$ the off-diagonal ones at $(1, 1,..., 1)$. This
gives $q$, $Q\in\Bbb N^+$ satisfying 
$$q<Q.$$
For example, in the cubic case, $p=1$, $P=\sum_{i=1}^b a_i^2$, $P'=1$ and $P_k=\sum_{i\neq k}a_i^2$ 
giving $q=2$ and $Q=4$. 

So at $a=(1, 1,  ..., 1)$, we have the following derivative matrix: 
$$\frac{\partial \Omega}{\partial a}(1, 1, ..., 1)=\pmatrix q&Q&Q&\cdots&Q\\Q&q&Q&\cdots&Q\\Q&Q&q&\cdots&Q\\\vdots&\vdots&\vdots&\vdots&\vdots\\Q&Q&Q&\cdots&q\endpmatrix.\tag 3.37$$
By inspection, the column vector with all entries $1$ is an eigenvector with eigenvalue $\lambda_1=q+(b-1)Q\neq 0$. Since the matrix has rank $2$, the 
other eigenvalue $\lambda_2$ is $(b-1)$- fold degenerate. Using the trace, we have 
$\lambda_1+(b-1)\lambda_2=bq$, so $\lambda_2=q-Q\neq 0$. Therefore
$$\det (\frac{\partial \Omega}{\partial a})(1, 1, ..., 1)\neq 0=\det (\frac{\partial \Omega}{\partial a}) (0, 0, ..., 0).$$
Hence $\det (\frac{\partial \Omega}{\partial a})$ is not a constant. 
Since $\det (\frac{\partial \Omega}{\partial a})$ is a polynomial in $a$ of degree at most
$(2p-1)b$, using the argument in (3.16) and  $\Vert \frac{\partial \Omega}{\partial a}\Vert=\Cal O(1)$ from (3.36, 3.32), this proves 
(3.33, 3.34) similar to the proof of (3.2), after a $\Cal O(\delta^{2p+2})$ perturbation. (Cf. the frequency 
modulation formula and an other proof of diffeomorphism in \cite{PP1}, Propositions 4 and 5.) 

Finally using the vector field 
$$\iota (\omega)=\sum_{k=1}^b\frac{n_k}{{\sum_{k=1}^b n_k^2}}\frac{\partial} {\partial \omega_k},$$
as well as (3.33), proves (3.35).
\hfill$\square$
\enddemo

\bigskip 
\head{\bf 4. The second step}\endhead
After the first step, the frequency $\omega$ is set at $\omega=\omega^{(1)}$. The 
modulation $\Delta \omega^{(1)}=\Cal O(\delta^{2p})$ is the same order 
as the matrix $A_0$ defined in (1.15). So it is not yet in a non-resonant form. Moreover the matrix 
$PA_0P$ defined by (1.17, 1.15) loses the restricted ``convolution'' structure, since the diagonal elements,
$$\aligned &\pm n\cdot\omega^{(1)}+j^2+(A_0)_{11}(0, 0)\\
=&\pm n\cdot\Delta \omega^{(1)}+(A_0)_{11}(0, 0),\endaligned$$
{\it depend} on $n$. (Recall that $(A_0)_{11}(0, 0)$ denotes the diagonal element of $(A_0)_{11}$.) 

However, since $\Delta \omega^{(1)}=\Cal O(\delta^{2p})$, in a region where $n$ is such that 
$|n\cdot \Delta \omega^{(1)}|\ll 1$, the resonance structure remains the same and there is still  
the block diagonal structure exhibited in Lemma 3.2. The goal of this section is to transform the 
$P$-equations into a non-resonant system. This is achieved in Lemma 4.1. The  block matrices ``near'' the origin
are dealt with using their determinants; while the others, variational arguments using the modulated 
frequency $\omega^{(1)}$.  Lemma 4.2 is a technical lemma ensuring that the determinants of the
matrices near the origin are not identically zero under appropriate conditions. Proposition 4.3 is
the culminating result having transformed the NLS into an amenable non-resonant system.   

\subheading {4.1 Invertibility of the linearized operator}

Toward that end, define the truncated linearized operator $F'_N({\tilde u}^{(0)},{\tilde v}^{(0)})$ evaluated at $\omega^{(1)}$ as
 $$
\cases
F'_N({\tilde u}^{(0)},{\tilde v}^{(0)}) (x,y) &=F'({\tilde u}^{(0)},{\tilde v}^{(0)})(x,y),\quad \Vert x\Vert_\infty\leq N,\, \Vert y\Vert_\infty\leq N ,\\
&=0,\qquad\qquad\qquad\qquad\text{ otherwise},
\endcases\tag 4.1
$$
where $${\tilde u}^{(0)}=\sum_{k=1}^b a_k e^{-ie_k\cdot\omega^{(1)} t }e^{ij_k\cdot x}$$
has the {\it modulated} frequency $\omega^{(1)}$; while the Fourier supports: 
$$\text{supp } {\tilde u}^{(0)}=\text{supp } u^{(0)}.$$
From Proposition 3.3, $$u^{(1)}={\tilde u}^{(0)}+\Cal O(\delta^{3-\epsilon}).$$
We note that $F'_N$ is $F'$ restricted to the set 
$$[-N, N]^{b+d}\cup[-N, N]^{b+d}\sim  [-N, N]^{b+d}\times\Bbb Z_2.$$ 

There is the analogue of Lemma 3.1.
\proclaim{Lemma 4.1} Let $\Cal B_{\epsilon, \delta}$ be the set defined (constructed)  in Proposition 3.3. Assume that  ($\flat$, $\flat\flat$) (renamed as (4.6, 4.8) below) hold. 
Let $N=|\log \delta|^s$ ($s>1$). There exists $\delta_0\in(0, 1)$,
such that  for all $\delta\in (0, \delta_0)$, there is a set 
$\tilde\Cal B_{\epsilon, \delta}$, $(0, 1]^b=\Cal B\supset\tilde\Cal B_{\epsilon, \delta}\supset\Cal B_{\epsilon, \delta}$, with 
$$\text{meas }\tilde\Cal B_{\epsilon, \delta}<\epsilon'.$$
If $a\in\Cal B\backslash\tilde \Cal B_{\epsilon, \delta}$,
then 
$$\Vert [F'_N({\tilde u}^{(0)}, {\tilde v}^{(0)})]^{-1}\Vert\leq \Cal O(\delta^{-2p-\epsilon})\tag 4.2$$
and there exists $\beta\in (0,1)$ such that
$$|[F'_N({\tilde u}^{(0)}, {\tilde v}^{(0)})]^{-1}(x,y)|\leq \delta^{\beta|x-y|}=e^{-\beta|\log\delta||x-y|}\tag 4.3$$
for all $|x-y|>1/\beta^2$.
\endproclaim 


\noindent{\it Remark.}  We could have skirted Lemma 3.1, by first solving the $Q$-equations in (1.13), cf. (3.36),
and directly arrived at Lemma 4.1. (This is done, e.g., in \cite{W5}.)
But since the proof of Lemma 4.1 is more technical, and this is the first paper using the method, Lemma 3.1 
is retained in order to illustrate some of the main issues. 

To prove Lemma 4.1, we note that since  
the Fourier support of $\tilde u^{(0)}$ satisfies
$$\text{supp }\tilde u^{(0)}=\text{supp }u^{(0)},$$
$$PF'_N(\tilde u^{(0)})P=\delta^{2p}\oplus_{\alpha\subset\Cal C}\Gamma_{\alpha}(\tilde u^{(0)})+\Cal O(\delta^{2p+2}),\tag 4.4$$ 
where $P$ as in (3.5), $F_N'$ as in (4.1), the decomposition is the same as in (3.7) with $\alpha$ connected sets on $\Cal C$ of sizes at most $2b+2d$;
moreover {\it as matrices}
$$\aligned \Gamma_\alpha(\tilde u^{(0)})&=\pmatrix \text{diag }(n\cdot\Omega)&0\\0&\text{diag }(-n\cdot\Omega)\endpmatrix+A'_\alpha(\tilde u^{(0)}, \omega^{(1)})\\
&=\pmatrix \text{diag }(n\cdot\Omega)&0\\0&\text{diag }(-n\cdot\Omega)\endpmatrix+A'_\alpha(u^{(0)}, \omega^{(0)}),\endaligned \tag 4.5$$
(note that the second $A'_\alpha$ is evaluated at $(u^{(0)}, \omega^{(0)}$)),
where $\Omega=\{\Omega_k\}$ as defined in (3.36) and $n$ are such that $(n, j)\in\alpha$. 
\smallskip
\noindent{\it The matrices near the origin}

We begin by studying the matrices near the origin. To prove
their invertibility, one starts by showing that the determinant polynomials are not identically zero
by showing that they are diagonally dominant at $a= (1, 1, ..., 1)$ under appropriate conditions. 
We prove

\proclaim{Lemma 4.2} Assume that $u^{(0)}$ of $b$ frequencies is generic satisfying the genericity conditions (Gi-iii),
and that 
$$n\cdot\omega^{(0)}\neq 0, \tag 4.6$$
for  $n\in\Bbb Z^b$, $0\neq |n|\leq N'=N'(p, d, b, \epsilon')$,
and $N'$ 
is independent of $\delta$ and assumed to be large. 
Then for all $\alpha$ in (4.4, 4.5), such that $$\alpha\subset\{ [ -N'/2, N'/2]^b\times [-N, N]^d\}\times \Bbb Z_2\backslash\Cal S,\tag 4.7$$ 
where $\Cal S$ as defined in (1.11), $\alpha$ satisfies 
$$|\alpha|\leq 2d+2.$$
Under the additional assumption 
$$b>C_pd\tag 4.8$$ for some $C_p>1$,
$$\det\Gamma_\alpha\not\equiv 0.$$
\endproclaim
\demo{Proof} 
We first show that under the condition (4.6),  for a connected set $\alpha$ satisfying (4.7), 
$$\pi\alpha\not\subseteq \{j_k\}_{k=1}^b\cup   \{-j_k\}_{k=1}^b,$$
where $\pi$ is the projection from $\Bbb Z^{b+d}$ to $\Bbb Z^d$. Therefore it follows from 
Lemma 3.2 that 
$$|\alpha|\leq 2d+2.\tag 4.9$$
This is derived by using contradiction. From Lemma 3.2,  if 
$$\pi\alpha\subseteq \{j_k\}_{k=1}^b\cup   \{-j_k\}_{k=1}^b,$$
then $$\pi\alpha= \{j_k\}_{k=1}^b\cup   \{-j_k\}_{k=1}^b.$$
Moreover, they are all of the form $E_\mu$ defined in ($\diamondsuit$),  
with $$\mu\cdot\omega^{(0)}=0.\tag 4.10$$ 
Since $\mu\neq 0$ for $\alpha$ satisfying (4.7), (4.10) contradicts (4.6). So (4.9) is satisfied for $\alpha$ in the 
decomposition (4.4), using the matrix equivalence
in (4.5).

Next we prove $$\det\Gamma_\alpha\not\equiv  0,$$by setting 
$$\quad a=\{a_k\}_{k=1}^b=(1, 1, ..., 1)$$
and showing that 
$$\det\Gamma_\alpha (1, 1, ..., 1)\neq 0,\tag 4.11$$
using (4.8, 4.9). It essentially consists of using the binomial formula to 
show that the matrix is diagonally dominant and an application of the Schur's  lemma.

From the first equality in (4.5) and the $Q$-equations in (3.36), the diagonal elements of $\Gamma_\alpha$ are
$$\align \Gamma_\alpha(n, j, n, j)&=\pm\sum_{k=1}^bn_k\Omega_k+(p+1)\text{ diag } (u*v)^{*p}\\
:&=(\Cal N+p+1)(u*v)^{*p} (-e_1, j_1, -e_1, j_1)+\Cal N(b-1)(u*v)^{*p}(-e_1, j_1, -e_2, j_2), \tag 4.12 \\
:&=(\Cal N+p+1)M_{11}+\Cal N(b-1)M_{12},\\
:&=D\endalign$$ 
where $u:=\tilde u^{(0)}$, $\Cal N=\pm\sum_{k=1}^bn_k$ and $n_k$ is the $k$th component of $n\in\Bbb Z^b$, with the plus sign for  $(n, j)\in\Cal C^+$, minus sign for  $(n, j)\in\Cal C^-$, and we used that $$a_1=a_2=\cdots=a_b,$$
to reach the second equality.

The off-diagonal elements of $\Gamma_\alpha$ are among the off-diagonals of 
$(p+1)(\hat u*\hat v)^{*p}$, $p \hat u*\hat u*(\hat u*\hat v)^{*(p-1)}$ and $p\hat v*\hat v*(\hat u*\hat v)^{*(p-1)}$, where
we have put back the hat with $\hat u$ the Fourier transform of $u$ and $\hat v$ the Fourier transform of $\bar u$. They are the 
Fourier coefficients of 
$(p+1)|u|^{2p}$, $pu^2|u|^{2(p-1)}$ and $pv^2|u|^{2(p-1)}$. Therefore it suffices to study the ``polynomials'' (product Fourier series) 
$$P_1=(|a_1|^2+|a_2|^2+....+|a_b|^2+\sum_{k\neq\ell} \frak a_k\bar {\frak a}_\ell)^p,$$
and 
$$P_2=(\frak a_1^2+\frak a_2^2+....+\frak a_b^2+2\sum_{k>\ell} \frak a_k\frak a_\ell)(|a_1|^2+|a_2|^2+....+|a_b|^2+\sum_{k\neq\ell} \frak a_k\bar {\frak a}_\ell)^{p-1},$$
where $$\frak a_k\bar {\frak a}_\ell:=a_k\bar a_\ell e^{i(e_\ell-e_k)\cdot\omega^{(1)} t}e^{i(j_k-j_\ell)\cdot x},$$
$$\frak a_k\frak a_\ell:=a_ka_\ell  e^{-i(e_\ell+e_k)\cdot\omega^{(1)} t}e^{i(j_k+j_\ell)\cdot x},$$
and $e_k$ and $e_\ell$ are the unit vectors in $\Bbb Z^b$ in the $k$th and $\ell$th directions as before. 

Using the binomial formula, we have 
$$P_1=\sum_{m=0}^{p}(|a_1|^2+|a_2|^2+....+|a_b|^2)^{p-m}C^m_p(\sum_{k\neq\ell} \frak a_k\bar {\frak a}_\ell)^m.\tag 4.13$$
The term that corresponds to $M_{11}$ is the Fourier coefficient $\hat P_1 (0,0)$ and to $M_{12}$, $\hat P_1(-e_1+e_2, j_1- j_2)$. 

To compute $M_{11}$, we expand the $m$-fold product. Assume $b\gg p$ is large. The combinatorial factor in front of the $\Cal O(b^p)$ term is 
$$1+2!C_p^2+3!C_p^3+...+m!C_p^m+ ... + p!C_p^p.$$
This is because in each term in the $m$-fold product, the indices must appear in pairs and are of the form: 
$$\frak a_\ell \bar {\frak a}_{k'} ... \frak a_k\bar {\frak a}_\ell,$$ 
(here we have made it  explicit on the index $\ell$).

Summing over the indices gives $\Cal O(b^m)$, the extra factor 
$m!$ comes from the number of ways of choosing the $m$ quadratic terms in the above product.
For example, when $m=2$, the factor $2$ in front of $C_p^2$ comes from 
the $2$ possible choices of $a_\ell \bar a_k$ ($k=k'$ here).

Similarly, the factor in front of the leading order $\Cal O(b^{p-1})$ term for $M_{12}$ is 
$$C_p^1+2!C_p^2+3!C_p^3+...+m!C_p^m+...+ p!C_p^p.$$
(Here two of the indices are fixed at $1$ and $2$ respectively, the other indices come in pairs as in $M_{11}$ and are summed over.)
So
$$\aligned D=&(\Cal N+p+1) (1+2!C_p^2+3!C_p^3+...+m!C_p^m+...+p!C_p^p)b^p\\
&+\Cal N(C_p^1+2!C_p^2+3!C_p^3+...+m!C_p^m+...+p!C_p^p)b^p+\Cal O(b^{p-1}).\endaligned$$
Setting the $\Cal O(b^p)$ term in $D$ to be $0$ gives 
for $p=1$, $\Cal N=-1$ and for $p\geq 2$,
$$\aligned \Cal N&=-\frac{(p+1)(\sum_{m=2}^p m!C_p^m+1)}{(p+1)+2\sum_{m=2}^p m!C_p^m}\\
&=-\big(\frac{p+1}{2}\big)\big(\frac{1+\frac{1}{A}}{1+\frac{p+1}{2A}}\big),\\
&=-\big(\frac{p+1}{2}\big)-\big(\frac{p+1}{2}\big)\big[\big(1+\frac{1}{A}\big)\sum_{n=1}^\infty(-x)^n+\frac{1}{A}\big]\endaligned\tag 4.14$$
where $A=\sum_{m=2}^p m!C_p^m$ and $x=\frac{p+1}{2A}$. 

Since 
$$0<\big|\big(\frac{p+1}{2}\big)\big[\big(1+\frac{1}{A}\big)\sum_{n=1}^\infty(-x)^n+\frac{1}{A}\big]\big|<1/2 $$
from direct computation, 
$\Cal N\notin\Bbb Z$ for $p\geq 2$.
Therefore with the exception of $p=1$, the diagonal elements of $\Gamma_\alpha=\Cal O(b^p)$. 

To proceed, we note that, using similar arguments as for $M_{12}$, any off-diagonal term of $\Gamma_\alpha$ that is a Fourier coefficient of $P_1$, 
$\hat P_1(n, j)$, $(n, j)\in\Bbb Z^{b+d}\backslash\{(0, 0)\}$, 
can be bounded by 
$$\sum_{m=1}^p b^{p-m}C^m_p m! b^{m-1}\leq C'_pb^{p-1}\tag 4.15$$
for some $C'_p>1$. 
This is because $n$ specifies the number of fixed indices and the remaining ones come in 
pairs and are summed over. The minimum number of fixed indices is two as in $M_{12}$. 
Similarly the bound (4.15) holds for $\hat P_2$. For $p\neq 1$, since $|D|=\Cal O(b^p)$, 
and $\Gamma_\alpha$ is at most a $(2d+2)\times (2d+2)$ matrix, {\it independent} of $b$,
(4.15) proves that
for $$b>C_pd,$$
$$\Vert \Gamma_\alpha^{-1}(1, 1, ..., 1)\Vert \lesssim \Cal O(b^{-p})$$
by using Schur's lemma and (4.9). So 
$$\det \Gamma_\alpha (1, 1, ..., 1)\neq 0.$$
Hence 
$$\det \Gamma_\alpha\not\equiv 0.$$

For the exceptional cubic case, $p=1$, we only need to verify that when $\Cal N=-1$, the matrix $\Gamma_\alpha$ is invertible.
In that case, from (4.12), it has 
$1$ on diagonal. The off-diagonals are among the Fourier coefficients of $2|u|^2$ or $u^2$ or $v^2$. The Fourier coefficients of $2|u|^2$ are manifestly all even.
The only odd coefficients in $u^2$, $v^2$ come from the terms $\frak a_\ell^2$, ${\bar \frak a_\ell}^2$, $\ell=1, 2, ..., b$. But if these terms appear in $PA_0P$, then 
the following two equations must be satisfied for some $\ell\in \{1, 2, ..., b\}$ and $(n, j)\in\Bbb Z^{b+d}$:     
$$
\cases
(n\cdot\omega^{(0)}+j^2)=0,\\
-(n+2e_\ell)\cdot\omega^{(0)}+(j-2j_\ell)^2=0.
\endcases
$$
Adding the two equations leads to 
$$(j-j_\ell)^2=0.$$
So $j=j_\ell$. Therefore, one may write 
$n=-e_\ell+\mu$, $\mu\neq 0$, satisfying $\mu\cdot\omega^{(0)}=0$. 
But this is a contradiction for connected sets $\alpha$ satisfying  
(4.7), using (4.6). So all the off-diagonals of $\Gamma_\alpha$ are even, equal to 
$2$ or $0$.  Invertibility, for all $b$,  follows from the determinant 
formula.  ($\Gamma_\alpha$ is an identity matrix on $\Bbb Z/2\Bbb Z$.) This completes the proof
for all $p$.
\hfill$\square$
\enddemo

\noindent {\it Remark.} 
The condition (4.8) is only used in Lemma 4.2. When $p=1$, since $\det\Gamma_\alpha\not\equiv 0$ for all $b$, this implies,  
in particular, that for the cubic NLS, the Theorem holds for any number of frequencies $b$.

\bigskip
\demo {Proof of Lemma 4.1}
We identify the set of connected sets $\{\alpha\}$ on $$\Cal C\cap [-N, N]^{b+d}\times \Bbb Z_2$$ 
with the set $\{1, 2,..., K_1\}$, where $K_1=K_1(N)$. So
$$PF'_{N}P=\delta^{2p}\oplus_k\Gamma_k(a)+\Cal O(\delta^{2p+2}),\quad k\leq K_1(N),$$
where each $\Gamma_k$ is of the form in (4.5).

With the addition of the diagonal term, $\Gamma_k$ is no longer a convolution matrix, as mentioned earlier. Moreover 
$|n|\leq N=|\log\delta|^s$ ($s>1$) depends on $\delta$. So we need to proceed differently because of 
uniformity considerations in estimates of type (3.16). 

Fix $$N_0=C_0p (2b+2d),\, N'=3N_0,\tag 4.16$$
for some large $C_0=C_0(p, b, d, \epsilon')>0$ to be determined by (4.17).   
For a given $\Gamma_k$, define the support of $\Gamma_k$ to be 
$$\Bbb Z^{b+d}\times \Bbb Z^{b+d}\supset\text{supp }\Gamma_k=\{(\pi x,\pi y) |\Gamma_k(x,y)\neq 0\},$$
where $\pi$ is the projection onto $\Bbb Z^{b+d}$.

For matrices $\Gamma_k$, such that $$\text{supp }\Gamma_k\cap \{[-N_0, N_0]^{b+d}\times [-N_0, N_0]^{b+d}\}\neq\emptyset,$$
we proceed as in the proof of Lemma 3.1. There are at most $K_0$ (independent of $\delta$) of these matrices.
Let  $P_k=P_k(a)=\det \Gamma_k(a)$, using Lemma 4.2 with the $N'$ in (4.16), we have that there exist $C$, $c>0$, such that for all $0<\epsilon<1$,
$$\text{ meas }\{a\in \Cal B||P_k|<\delta^{\epsilon}, \text{ all }k\leq K_0\}\leq C\delta^{c\epsilon}.$$
So $\Vert \Gamma_k^{-1}\Vert\leq\Cal O(\delta^{-\epsilon})$ for all $k\leq K_0$. 

For matrices $\Gamma_k$ with $k>K_0$, 
$$\text{supp }\Gamma_k\cap \{[-N_0, N_0]^{b+d}\times [-N_0, N_0]^{b+d}\}=\emptyset$$
by definition. We use perturbation theory. For any $\Gamma_k$, fix $N''$, with $|N''|>N_0$, such that
for all 
$(n, j)\in\text{supp }\Gamma_k$, we can write
$(n, j)=(N'', 0)+(n'', j)$ with $|n''|\leq 2p(2b+2d)$.

Define the directional derivative $\frac{d}{d\Omega}$, where $\Omega$ is defined in (3.36),  to be 
$$\frac{d}{d\Omega}:=\sum_{i=1}^{b}\frac{N_i''}{\Vert N''\Vert_2^2}\cdot\frac{\partial}{\partial\Omega_i}.$$ 
Assume that $\Gamma_k$ is an $m\times m$ matrix, $m\leq 2b+2d$ from Lemma 3.2. The determinant $P_k$
can be written as:
$$P_k=(-1)^{m'} (N''\cdot \Omega)^m+q_{m-1} (N''\cdot\Omega)^{m-1}+...+q_0,$$
where $m'$ is the number of connected sites on $\Cal C^-$, $q_\mu=q_\mu(n'', a,\Omega)$, $\mu=0, 1, ..., m$,
and are {\it independent} of $N''$. 

Taking the $m$-th order derivative yields 
$$|\frac{d^m}{d\Omega^m}P_k|>\frac{1}{2}\tag 4.17$$
for $C_0$ large enough depending only on $p$, $b$, $d$ and $\epsilon'$, where we used the form of the matrix in (4.5) and that there are only finite 
types of  ``convolution'' matrices $A'_k$ and that  
$$\Vert \frac{\partial \Omega}{\partial a}\Vert\asymp \Vert\big( \frac{\partial \Omega}{\partial a}\big)^{-1}\Vert\asymp \Cal O_{\epsilon'}(1).$$
cf.  (3.32, 3.33). Since $\Vert \Gamma_k\Vert \leq \Cal O(|\log\delta|^s)$, the $m$-th variation at $P_k=0$ gives $\Vert \Gamma_k^{-1}\Vert~\leq ~\Cal O(\delta^{-\epsilon})$
for all $K_0<k\leq K_1(N)$ away from a set in $a$ of measure less than $\delta^{\epsilon/4(b+d)}$, where we also
used $K_1(N)\leq\Cal O(|\log\delta|^{2(b+d)s})$, $s>1$. 
Using the above estimates on $\Vert \Gamma_k^{-1}\Vert$ in (4.4), the expression right above (3.14), with $F'$ replaced by $F'_N$,  then gives (4.2). 
Setting $\delta_0$ to satisfy $\delta_0^{\epsilon/4(b+d)}<\epsilon'/2$,  yields the measure estimate. This is as in the proof of Lemma 3.1.

Lastly, since the geometry of 
the resonant structure remains the same, the point-wise estimates on 
$[F'_N({\tilde u}^{(0)},  {\tilde v}^{(0)})]^{-1}$ can be obtained as in the proof of (3.3), yielding (4.3).
\hfill$\square$
\enddemo

\subheading {4.2 The second iteration}

Let $N=|\log\delta|^s$ for some $s>1$ and $F'_N$
be $F'(\omega^{(1)}, {\tilde u}^ {(0)}, {\tilde v}^ {(0)})$ restricted to the set 
$$[-N, N]^{b+d}\times \Bbb Z_2,\tag 4.18$$
as before. Clearly the set in (4.18)
contains $$\text{supp } F(\omega^{(1)}, {\tilde u}^ {(0)}, {\tilde v}^ {(0)})$$
as a subset for small $\delta$. 

To prepare for the upcoming inductive construction, we redefine $\Delta u^{(1)}$ 
to be 
$$\aligned\pmatrix\Delta u^{(1)}\\\Delta v^{(1)}\endpmatrix=&[F_N']^{-1}(\omega^{(1)}, {\tilde u}^ {(0)}, {\tilde v}^ {(0)})F(\omega^{(1)}, {\tilde u}^ {(0)}, {\tilde v}^ {(0)})\\
:=&\pmatrix\Delta u^{(1)}(\omega^{(1)})\\\Delta v^{(1)}(\omega^{(1)})\endpmatrix,\endaligned\tag 4.19$$
which involves the same frequency (frequency at the same stage of iteration) and is more conducive to applying the implicit function theorem 
to the $Q$-equations. Clearly (3.29) remains valid after the redefinition with a possible lowering of $\beta$. (We note that the previous definition in (1.16) entailed  
$\Delta u^{(1)}=\Delta u^{(1)}(\omega^{(0)})$ instead.) As before (4.19) is defined on the domain of the $P$-equations, which is in the complement of the
set $\Cal S$ defined in (1.11). We keep the definition of $\omega^{(1)}$ in (1.18). 

In other words, the Newton scheme that we will use in sect. 5 unfolds as follows: $\omega^{(-1)}=\{j^2_k\}_{k=1}^b$,
$u^{(-1)}=0$;  $\omega^{(0)}=\omega^{(-1)}$ since $u^{(-1)}=0$, $u^{(0)}$ as in (1.8); $\omega^{(1)}$ as  in (1.18)
and $\Delta u^{(1)}=\Delta u^{(1)}(\omega^{(1)})$ as in (4.19) ...

We summarize the findings so far in the following amplitude-frequency 
modulation proposition. For simplicity, we use $u$ to denote both the function and its Fourier series as it should 
be clear from the context.
\proclaim
{Proposition 4.3}
Assume that 
$$u^{(0)}(t, x)=\sum_{k=1}^b  a_k e^{ij_k\cdot x}e^{-i{j^2_k}t},$$
a solution to the linear Schr\"odinger equation (1.2) is generic,
$a=\{a_k\}_{k=1}^b\in(0, 1]^b=\Cal B $. Let $\epsilon, \epsilon'\in (0, 1)$. 
There exists $\delta_0>0$, such that  for all $\delta\in (0, \delta_0)$, there is a set $\tilde\Cal B_{\epsilon, \delta}$, $(0, 1]^b=\Cal B\supset\tilde\Cal B_{\epsilon, \delta}\supset\Cal B_{\epsilon, \delta}$
(the set in Proposition 3.3), with 
$$\text{meas }\tilde\Cal B_{\epsilon, \delta}<\epsilon'.$$
There exists $\beta\in (0, 1)$ such that 
if $a\in\Cal B\backslash\tilde\Cal B_{\epsilon, \delta}$, then 
$$\Vert \Delta u^{(1)}\Vert_{\ell^2(\rho)} =\Vert \Delta v^{(1)}\Vert_{\ell^2(\rho)}= \Cal O(\delta^{3-\epsilon}),\tag 4.20$$
where $\Delta u^{(1)}$, $\Delta v^{(1)}$  as defined in (4.19) and $\rho$ is a weight on $\Bbb Z^{b+d}$ satisfying 
$$\aligned \rho(x)=&e^{\beta|\log\delta| |x|}, \text{ for } |x|>1/\beta^2,\\ 
=&1,\qquad\qquad\text{  for } |x|\leq 1/\beta^2.\endaligned$$
$$\aligned&\Vert \Delta \omega^{(1)}\Vert=\Vert \omega^{(1)}-\omega^{(0)}\Vert\asymp \delta^{2p},\\
&\Vert \frac{\partial \omega^{(1)}}{\partial a}\Vert\asymp \delta^{2p},\\
&\Vert(\frac{\partial \omega^{(1)}}{\partial a})^{-1}\Vert\lesssim \Cal O_{\epsilon'}(\delta^{-2p}),\\
&\big| \det(\frac{\partial \omega^{(1)}}{\partial a})\big|\gtrsim \Cal O_{\epsilon'}(\delta^{2pb}),\endaligned$$
(as in Proposition 3.3, $\frac{\partial \omega^{(1)}}{\partial a}$ is meant in the classical sense),
and $\omega^{(1)}$ is Diophantine
$$\Vert n\cdot \omega^{(1)}\Vert_{\Bbb T}\geq\frac{\kappa\delta^{2p}}{|n|^\gamma},\quad n\in\Bbb Z^b\backslash\{0\},\,\kappa>0, \gamma>2b+1,$$
where $\Vert\,\Vert_{\Bbb T}$ denotes the distance to integers in $\Bbb R$, $\kappa$ and 
$\gamma$ are independent of $\delta$.

We have moreover,
$$\Vert F(\omega^{(1)}, u^{(1)}, v^{(1)})\Vert_{\ell^2(\rho)\times\ell^2(\rho)}=\Cal O(\delta^{2p+5-2\epsilon}),\tag 4.21$$
$$\Vert [F'_N(\omega^{(1)}, u^{(1)}, v^{(1)})]^{-1}\Vert\leq \Cal O(\delta^{-2p-\epsilon}),\tag 4.22$$
and 
$$|[F'_N(\omega^{(1)}, u^{(1)}, v^{(1)})]^{-1}(x,y)|\leq \delta^{\beta|x-y|}=e^{-\beta|\log\delta||x-y|}\tag 4.23$$
for all $|x-y|>1/\beta^2$.
\endproclaim 
\demo{Proof} We only need to prove (4.20, 4.21), the rest reiterates Proposition 3.3 and Lemma 4.1.
The estimates in (4.20) follow from (4.2). To prove (4.21), we write 
$$\aligned F(u+\Delta u)&= F(u)+F'(u)\Delta u + \Cal O (\Vert F''(\bar u)\Vert \Vert \Delta u\Vert^2),\\ 
&=(F'-F'_N)[F'_N]^{-1}F(u)+\Cal O(\delta^{2p+5-2\epsilon})\\
&=\Cal O(\delta^{2p+5-2\epsilon}), \endaligned\tag 4.24$$
where $u$ stands for $\pmatrix {\tilde u}^{(0)}\\{\tilde v}^{(0)}\endpmatrix$,
$\Delta u$ stands for $\pmatrix \Delta u^{(1)}\\\Delta v^{(1)}\endpmatrix$
and we used (4.20, 4.22, 4.23).
\hfill $\square$ 
\enddemo
\bigskip 
\head{\bf 5. Proof of the Theorem}\endhead
Proposition 4.3 puts the construction in a non-resonant form with $\omega^{(1)}$ as the 
parameter.  It provides the input for the initial scale in the Newton scheme in \cite{B3}. To 
continue the iteration, we need the analogues of (4.22, 4.23) at larger scales. This is attained as follows.

Let $T=F'$ be the linearized operator defined as in (1.14-1.15) and the restricted operator 
$T_N=F'_N$ as defined in (4.1). To increase the scale from $N$ to a larger scale $N_1$, we pave the $N_1$
cubes with $N$ cubes. In the $j$ direction, this is taken care of by perturbation; while in the $n$
direction by adding an additional parameter $\theta\in\Bbb R$ and consider $T(\theta)$: 
$$
T(\theta) =\pmatrix \text {diag }(n\cdot\omega+j^2+\theta)&0\\ 0& \text {diag }(-n\cdot\omega+j^2-\theta)\endpmatrix+\delta^{2p} A,$$
where $\delta^{2p}A$ correspond to the $A$ defined in (1.15). (Recall the rescaling $a\to\delta a$ starting in sect. 3.) 

\noindent{\it Remark.} This one dimensional parameter $\theta$ is merely an {\it auxiliary} variable. Using the covariance
of $n\cdot\omega+\theta$, all estimates in $\theta$ will be transformed into estimates in $\omega$ in the
Newton scheme construction of $u$ and $\theta$ is always {\it fixed} at $0$ there. So in particular, $A=A(\omega, u, v)$ is a T\"oplitz matrix 
{\it independent} of $\theta$. 

\subheading {5.1 The $\theta$ estimates}
  
Let $N=|\log\delta|^s$ $(s>1)$ as in Proposition 4.3 and $T_N(\theta)=T_N(\theta; u^{(1)}, v^{(1)})$ evaluated at
$\omega^{(1)}$. We have the following estimates. 
\proclaim{Lemma 5.1} Assume that $u^{(0)}=\sum_{k=1}^b  a_k e^{ij_k\cdot x}e^{-i{j^2_k}t}$ 
a solution to the linear Schr\"odinger equation (1.2) is generic and $a\in\Cal B\backslash\tilde\Cal B_{\epsilon, \delta}$, the set defined in Proposition 4.3. Then
there exists $\delta_0>0$, such that for all $\delta\in(0, \delta_0)$
$$\Vert [T_N(\theta)]^{-1}\Vert\leq \Cal O(\delta^{-2p-\epsilon})<e^{N^\sigma}\tag 5.1$$ 
for some $\sigma\in (0, 1)$ 
and there exists $\beta\in (0,1)$ such that
$$|[T_N(\theta)]^{-1}(x,y)|\leq \delta^{\beta|x-y|}=e^{-\beta|\log\delta||x-y|}\tag 5.2$$
for all $|x-y|>1/\beta^2$, away from a set $B_N(\theta)\subset\Bbb R$ with 
$$\text{meas }B_N(\theta)<\delta^{2p+c\epsilon}<e^{-N^\tau}, \tag 5.3$$
for some $\tau \in (0, 1)$.
\endproclaim 
\demo{Proof} We assume $\omega=\omega^{(1)}$ is fixed.  Since $|n|\leq N=|\log\delta|^s$ $(s>1)$, the spectrum of $T_N$:
$$\sigma(T_N)\subset \bigcup \{\Theta+\delta^{2p}I\},$$
where $\Theta\in\Bbb Z$, 
$$I=[-C |\log\delta|^{s}, C |\log\delta|^{s}],$$ 
for some $C>0$, and the union is over $\Theta$ such that $|\Theta|\leq 2d |\log\delta|^{2s}$.
So it suffices to look at $\theta$ such that 
$$\theta\in  \bigcup \{\Theta+\delta^{2p}[-2C |\log\delta|^{s}, 2C |\log\delta|^{s}]\}.$$
Write $\theta=\Theta+\delta^{2p}\theta'$, then 
$$\aligned
T_N(\theta) =&\pmatrix \text {diag }(n\cdot\omega^{(0)}+j^2+\Theta)&0\\ 0& \text {diag }(-n\cdot\omega^{(0)}+j^2-\Theta)\endpmatrix\\
&+\delta^{2p} \pmatrix \text {diag }(n\cdot\tilde\omega+\theta')&0\\ 0& \text {diag }(-n\cdot\tilde\omega-\theta')\endpmatrix
+\delta^{2p} A_N,\endaligned$$
where $\omega^{(0)}\in\Bbb Z^b$, $\Theta\in\Bbb Z$, $\tilde\omega=\Delta\omega^{(1)}/\delta^{2p}$ is Diophantine
and $A_N$ is the restricted $A$ as defined in (1.15, 4.1).

Let $$\Cal H=\delta^{2p}\left[ \pmatrix \text {diag }(n\cdot\tilde\omega+\theta')&0\\ 0& \text {diag }(-n\cdot\tilde\omega-\theta')\endpmatrix
+A_N\right].$$
Let $P_+$ be the projection onto the set $\{(n,j) |n\cdot\omega^{(0)}+j^2+\Theta=0\}$ and 
$P_-$ the projection onto the set $\{(n,j) |-n\cdot\omega^{(0)}+j^2-\Theta=0\}$ when $\Theta\neq 0$;
when $\Theta=0$, use the definition in (3.5). Define 
$$P=\pmatrix P_+&0\\0&P_-\endpmatrix$$
and $P^c$ the projection onto the complement as before. 

We proceed using the Schur reduction as 
in the proofs of Lemmas 3.1 and 4.1. It suffices to estimate
$P\Cal H P=\oplus_k\Gamma_k(\theta)$
as $P^cT_NP^c$ is invertible, $\Vert (P^cT_NP^c-\lambda)^{-1}\Vert\leq 4$ uniformly in $\theta$ and
$\lambda\in[-1/4, 1/4]$.

Since $\det\Gamma_k$ is a polynomial in $\theta'$ of degree at most $(2b+2d)$ and the highest degree term 
has coefficient $\pm 1$, we obtain that 
$$\text{ meas }\{\theta'|\Vert\Gamma_k^{-1}(\theta')\Vert>\frac{1}{2}\delta^{-2p-\epsilon}\}\leq C\delta^{\frac{\epsilon}{4(b+d)}}\qquad (\epsilon>0),\tag 5.4$$
where we also used $\Vert \Gamma_k\Vert\leq\Cal O(|\log\delta|^s)$. Summing over $\Theta$ and the number of possible $\Gamma_k$ and taking into account the $\delta^{2p}$ factor in front of $\theta'$, 
we then obtain
$$\text{ meas }\{\theta\Vert T_N^{-1}(\theta)\Vert>\delta^{-2p-\epsilon}\}\leq C\delta^{2p+\frac{\epsilon}{8(b+d)}}\qquad (\epsilon>0),\tag 5.5$$
which gives (5.1, 5.3). The pointwise estimate (5.2) follows as in the proof of (3.3) in Lemma 3.1.
\hfill$\square$ 
\enddemo
\smallskip
Lemma 5.1 enables us to apply the inductive Lemma 19.38-19.65 in \cite {Chap 19, B3} to obtain the corresponding estimates at larger scales. 
The point-wise bounds in (19.15', 19.66), which follow from an application of
Lemma 7 in \cite{B1}, should be replaced by sub-exponential instead. This is because 
the sizes of the resonant clusters are, in general, much larger than their separations, cf. Lemma 19.10 in \cite{B3}.
(However, the proofs in \cite{B3} are not affected, cf. \cite{W1} for another instance where both the norm and point-wise bounds are 
sub-exponential.) 

\subheading {5.2 The Newton construction and proof of the Theorem} 

The proof of the Theorem is an induction. This is essentially the same as the combination of Chaps. 18 and 19 of \cite {B3}, cf. also \cite{B1}.
We first lay down the induction hypothesis.
Let $$a\in (0, 1]^b,\, \omega\in \delta^{2p}(-B, B)^b+(j^2_1, j^2_2, ..., j^2_b),\tag 5.6$$
where $B=B(p, b)$ and define
$$\tilde\omega=[\omega- (j^2_1, j^2_2, ..., j^2_b)]/\delta^{2p}\in (-B, B)^b.\tag 5.7$$

Let $M$, $R$ be large integers.  On the {\it entire} $(a,\tilde\omega)$ space, namely $(0, 1)^b\times (-B, B)^b$, assume that the following is
satisfied for $r\in [1, R]$: 
\item{(Hi)} $\text{supp }u^{(r)}\subseteq B(0, M^r)$ ($\text{supp }u^{(0)}\subset B(0, M)$),
\item{(Hii)} $\Vert \Delta u^{(r)}\Vert<\delta_r$,
$\Vert \partial \Delta u^{(r)}\Vert<\bar\delta_r$ with $\delta_{r+1}\ll\delta_{r}$ and $\bar\delta_{r+1}\ll\bar\delta_{r}$,

where $\partial$ refers to derivations in $a$ or $\tilde\omega$ and 
$\Vert\,\Vert:=\sup_{a,\tilde\omega}\Vert\,\Vert_{\ell^2(\Bbb Z^{b+d})\times \ell^2(\Bbb Z^{b+d})}$
\item{(Hiii)} $|u^{(r)}(\xi)|<e^{-|\xi|^c}$ for some $c\in (0, 1)$,
\bigskip
Using (Hi-iii), an application of the implicit function theorem to the
$Q$-equations: 
$$\tilde\omega_k(a)=\frac{[{(u*v)}^{*p}*u](-e_k, j_k)}{a_k}+\delta^2\frac{\hat H(-e_k, j_k)}{a_k}\tag 5.8$$
with $u=u^{(r)}$ yields
$$ {\tilde\omega}^{(r)}(a)=\Omega(u^{(0)}(a))+\delta^{2-\epsilon} \phi_r(a)+\delta^2 P(u^{(0)}(a)),\tag 5.9$$
where $0<\epsilon<1$ and $\Vert\partial\phi_r\Vert<C$. We define $\phi_0=0$ and denote the graph of $ {\tilde\omega}^{(r)}$ by $\Phi_r$. 
The vector valued polynomial $\Omega=\{\Omega_k(u^{(0)}(a))\}_{k=1}^b$ is as in (3.36) with each component homogeneous in $a$ of degree $2p$.
Similarly $P=\{P_k(u^{(0)}(a))\}_{k=1}^b$ with each component $P_k$ a polynomial in $a$ of bounded degree at least $2p+2$.
Moreover  by (Hii),
$$| {\tilde\omega}^{(r)}-  {\tilde\omega}^{(r-1)} |\lesssim \Vert u^{(r)}-u^{(r-1)}\Vert <\delta_r,\tag 5.10$$
so that $\Phi_{r-1}$ is a $\delta_r$ approximation of $\Phi_{r}$. This can be seen as follows.

Consider the right side of (5.8) as a function of $(a, \tilde\omega)$ and rewrite (5.8) as
$$\Cal F_k (a, \tilde\omega)=0,$$
for $k=1, ..., b$. Since $$u^{(r')}(a,\tilde\omega)=u^{(0)}(a)+\sum_{i=1}^{r'}\Delta u^{(i)}(a,\tilde\omega),$$
$\Cal F_k$ may be written in the form:
$$\Cal F_k(a, \tilde\omega)=f_k(a, \tilde\omega)+\Omega_k(u^{(0)}(a))+\delta^2 P_k(u^{(0)}(a)).$$
Let $X$ and $Y$ be the partial derivative matrices:
$$X=[[\frac{\partial {f_k}}{\partial a_\ell}]]\text{ and } Y=[[\frac{\partial {f_k}}{\partial\tilde\omega_\ell}]],\, k, \ell =1, ..., b.$$
The hypothesis (Hii) gives 
$$X=\Cal O(\sum_{i=1}^{r'}\Delta u^{(i)})+\Cal O( \sum_{i=1}^{r'}\partial_a\Delta u^{(i)})$$ and 
$$Y=\Bbb  I+\Cal O(\sum_{i=1}^{r'}\partial_{\tilde \omega}(\Delta u^{(i)})),$$
where the $\Cal O$ depends on $u^{(0)}$ and (or) $\partial_a u^{(0)}$.

So the partial derivative matrix 
$$\aligned \frac{\partial \tilde \omega}{\partial a}:&=[[\frac{\partial {\tilde\omega_k}}{\partial a_\ell}]], \, k, \ell =1, ..., b\\
&=-Y^{-1}(X+\partial \Omega/\partial a+\delta^2 \partial P/\partial a)\endaligned$$
is well-defined. The difference matrix satisfies 
$$\Vert  \frac{\partial {\tilde \omega}^{(r)}}{\partial a}- \frac{\partial {\tilde \omega}^{(r-1)}}{\partial a}\Vert \lesssim \Vert\partial _a\Delta u^{(r)} \Vert +\Vert\Delta u^{(r)} \Vert ,$$
since  $$\Delta\tilde \omega^{(r)}(0)=\Delta u^{(r)}(0)=0,$$
and $|a|\leq 1$,
this proves (5.10).
 \hfill $\square$
\bigskip
Below we continue with the assumptions on the {\it restricted} intervals in $(a,\tilde\omega)$ on $(0, 1)^b\times (-B, B)^b$, where one could construct approximate solutions.

\item{(Hiv)} There is a collection $\Lambda_r$ of intervals of size $cM^{-r^C}\delta^\epsilon$, $0<\epsilon<1$,  such that 
\item{(a)} On $I\in\Lambda_r$, $u^{(r)}(a,\tilde\omega)$ is given by a rational function in $(a,\tilde\omega)$ of degree at most 
$M^{Cr^3}$, for some $C>1$
\item{(b)} For $(a,\tilde\omega)\in\bigcup_{I\in\Lambda_r} I$,

$\Vert F(u^{(r)})\Vert<\kappa_r$,
$\Vert \partial F(u^{(r)})\Vert<\bar\kappa_r$ 
with $\kappa_{r+1}\ll\kappa_{r}$ and $\bar\kappa_{r+1}\ll\bar\kappa_{r}$
\item{(c)} Let $N=M^r$. For $(a,\tilde\omega)\in\bigcup_{I\in\Lambda_r} I$, $T=T(u^{(r-1)}):=F'(u^{(r-1)})$ satisfies

$\Vert T_N^{-1}\Vert <M^{(r^C+|\log\delta|)}$, 

$|T_N^{-1}(\xi, \xi')|<e^{-|\xi-\xi'|^c}$ for $|\xi-\xi'|>Cr^{C/c}$,

where $T_N$ is $T$ restricted to $$[-N, N]^{b+d}\cup [-N, N]^{b+d}\sim  [-N, N]^{b+d}\times \Bbb Z_2.$$
\item{(d)} Each $I\in\Lambda_r$ is contained in an interval $I'\in\Lambda_{r-1}$ and 
$$\text{meas}_b (\Phi_r\cap(\bigcup_{I'\in\Lambda_{r-1}} I'\backslash \bigcup_{I\in\Lambda_{r}}I)<\delta^{c\epsilon}[\exp\exp(\log (r+1))^{1/3}]^{-1},\, r\geq 2.$$
If $\tilde\omega\in \Phi_r\cap I$, then
$$\Vert n\cdot\tilde\omega\Vert_{\Bbb T}\geq \frac{\tilde\kappa}{|n|^\gamma},\, \tilde\kappa>0,\, \gamma>2b+1$$
for $|n|\leq M^r$ after identification of $\Phi_r\cap I$ with an interval in $\Bbb R^b$. 

\proclaim{Lemma 5.2} Assume that 
$$u^{(0)}(t, x)=\sum_{k=1}^b  a_k e^{ij_k\cdot x}e^{-i{j^2_k}t},$$
a solution to the linear Schr\"odinger equation (1.2) is generic. Let $R=|\log\delta|^c$ for some $c\in (0, 1)$. Then the induction hypothesis (Hi-iv) are satisfied for $r\in [1, R]$ and small $\delta$
with $\delta_r+\bar\delta_r$ and $\kappa_r+\bar\kappa_r$
satisfying $$\aligned \log\log\frac{1}{\delta_r+\bar\delta_r}&\sim r,\\
\log\log\frac{1}{\kappa_r+\bar\kappa_r}&\sim r.\endaligned\tag 5.11$$
\endproclaim

\demo{Proof} Since $N=M^R=M^{|\log\delta|^c}\ll\delta^{-1}$ for $0<c<1$, $$N\cdot\Delta\omega^{(R)}\sim N\cdot\Delta\omega^{(1)}\ll 1,$$ the resonance structure remains the same.  The Lemma
follows by repeating the construction in Proposition 4.3 (Lemma 4.1) $R$ times using the modified Newton scheme as in (4.24), which still leads to double exponential convergence because of the 
point-wise exponential estimates on $[F'_N]^{-1}$. The extension to the entire $(a,\tilde\omega)$ space is done as in \cite{sect. 10, B2}, in particular,  (10.33-10.37).   

Below we amplify the derivative estimates in (Hii, iv, b), the estimate on the size of the intervals in (Hiv), the bound on the degree of the rational functions in (Hiv, a) and the measure estimates in (Hiv, d), as the rest 
are direct products of the Newton construction. 

For the derivative estimates in (Hii, iv, b), we use the formula 
$$\Delta u^{(r)}=-[F'_N(u^{(r-1)})]^{-1}F(u^{(r-1)})\tag 5.12$$
and assume (Hii, iv, b) and (5.11) are satisfied at stage $r-1$. Taking the derivatives, we then have
$$\aligned \partial \Delta u^{(r)}=&-[F'_N(u^{(r-1)})]^{-1}\partial F(u^{(r-1)})\\
&\quad +[F'_N(u^{(r-1)})]^{-1}\partial [F'_N(u^{(r-1)})][F'_N(u^{(r-1)})]^{-1}F(u^{(r-1)})\\
=&\Cal O(M^{r^C}\cdot\bar\kappa_{r-1})+\Cal O(M^{2r^C}\cdot\kappa_{r-1})\\
:=&\bar\delta_r,\endaligned$$
which gives the derivative estimates in (Hii) at stage $r$. Using (4.24) with $u=u^{(r-1)}$ and $\Delta u= \Delta u^{(r)}$ and also taking the
derivatives of (4.24) and using the second bound in (5.11) at stage $r-1$ then give:
$$\aligned \partial F(u^{(r)})=&\{\partial (F'-F'_N)[F'_N]^{-1}F\}(u^{(r-1)})\\
&+\{(F'-F'_N)[F'_N]^{-1}\partial F'_N [F'_N]^{-1}F\}(u^{(r-1)})\\
&+\{(F'-F'_N)[F'_N]^{-1}\partial F\}(u^{(r-1)})\\
&+\Cal O(\partial F''(\Delta u^{(r-1)})^2)+\Cal O(F''(\Delta u^{(r-1)})(\partial\Delta u^{(r-1)}))\\
=&\Cal O(e^{-M^r}M^{2r^C}\bar\kappa_{r-1})+\Cal O(\delta_{r-1}\bar \delta_{r-1})\\
:=&\bar\kappa_r.\endaligned$$
This yields (Hiv, b) at stage $r$. So (5.11) and moreover the bounds in Lemma 5.5 are satisfied for $r\leq R$. The $Q$-equations express the modulated frequencies $\omega^{(r)}$ in a $\delta$ series for 
$r\leq R$.

The condition on the size of the intervals in (Hiv) is satisfied by using the stability of the estimates in (Hiv, c) under perturbations of the same size.
The measure estimates in (Hiv, d) follow from that for $2\leq r\leq R=|\log\delta|^c$, $0<c<1$, the additional excision verifies
$$ \text{meas }\tilde\Cal B_{\epsilon, \delta}\backslash \Cal B_{\epsilon, \delta}=\delta^{\tilde c\epsilon}\ll \delta^{c\epsilon}\exp\exp(\log\log\delta)^{1/3}$$ for $1>\tilde c>c>0$ and that the map $a\mapsto\omega(a)$ is diffeomorphic, cf. 
Lemma 3.1 and Proposition 4.3. 

From (5.12), $u^{(r)}$ is a rational function in $(a,\tilde\omega)$. Expressing the matrix elements of $[F'_N]^{-1}$ with $N=M^r$ as a ratio of determinants, which are polynomials in $(a,\tilde\omega)$ gives the bound on the degree in (Hiv, a) and concludes the proof. 
\hfill$\square$
\enddemo

Let $u$ denote $u^{(0)}$, $u^{(1)}$, ... For all $\bar N$, let $T_{\bar N}=T_{\bar N}(u)$ be the linearized operator 
evaluated at $u$ and restricted to $\{j+[-\bar N, \bar N]^{b+d}\}\times \Bbb Z_2$, where $j\in\Bbb Z^d$. (For simplicity the 
$j$ subindex is omitted.) Define the operator $T_{\bar N}(\theta)$ as before. Assume that (Hi-iv) hold at stage $r$. 
When $|j|\leq 2\bar N$, on the set of intervals $\Lambda_r$ in (Hiv), there is moreover the following estimates.

\proclaim{Lemma 5.3} There exist $c$, $\sigma$, $\tau\in (0, 1)$ ($c>\sigma>\tau$) such that 
$$\aligned \Vert T_{\bar N}^{-1}(\theta)\Vert &<e^{{\bar N}^\sigma},\\
|T_{\bar N}^{-1}(\theta)(x, y)| &< e^{-|x-y|^c}\endaligned\tag 5.13$$
for all $|x-y|>\bar N/10$, away from a set $B_{\bar N}(\theta)$ with 
$$\text{meas } B_{\bar N}(\theta)<e^{-\bar N^\tau},$$
where $u=u^{(r)}$, $|\log\delta|^s\leq \bar N\leq r^C$, $s$, $C>1$, $r\geq R$. In fact (5.13) holds for all $u= u^{(r')}$ with $r'>r\geq R$ fulfilling 
assumptions (Hi-iii) and verifying (5.11).
\endproclaim

\demo{Proof}
Set $N_0=|\log\delta|^s$, $s>1$ and $N_1=e^{N_0^{c'}}=e^{|\log\delta|^{sc'}}$ with $c'\in(0, 1)$ such that $0<sc'<c<1$, the same $c$ as in the definition of $R$ in Lemma 5.2. So the resonance
structure remains the same with $N_1\cdot\Delta\omega^{(R)}\ll 1$. Repeating the arguments in Lemma~5.1 with $u^{(1)}$ replaced by $u^{(R)}$ and for different intervals in $\Bbb Z^{b+d}\times \Bbb Z_2$, we then obtain (5.13) for the scales $\bar N\in[N_0, N_1]$.

Clearly at scales $\bar N\leq N_1$, we may replace
$u^{(R)}$ by any $u^{(r)}$ for $r>R$ if 
$$\Vert u^{(r)}-u^{(R)}\Vert \leq \Cal O(e^{-c''\bar N})$$
for some $c''>0$ and (5.13) remains valid.
From  (Hi-iv) and (5.11), the above bound is verified for $c'$ such that $0<sc'<c$. After possibly lowering the earlier $c'$ in the definition of $N_1$
and using the estimates for the scales in $[N_0, N_1]$ as the initial input, 
the induction lemmas in \cite {B3}, Lemma 19.38-65 and Lemma 19.13 then 
conclude the proof of (5.13) for all $r$. (Here one may assume that the $\Lambda_r$ in (Hiv) is constructed with the additional excision so that the bounds in Lemma 19.38-65 are available.)   
\hfill$\square$
\enddemo

There are related estimates on the set of intervals $\Lambda_r$ when $|j|>2\bar N$.  

\proclaim{Lemma 5.4} There exist $c$, $\sigma$, $\tau\in (0, 1)$ ($c>\sigma>\tau$) such that 
$$\aligned \Vert T_{\bar N}^{-1}(\theta)\Vert &<e^{{\bar N}^\sigma},\\
|T_{\bar N}^{-1}(\theta)(x, y)| &< e^{-|x-y|^c}\endaligned$$
for all $|x-y|>\bar N/10$, 
provided 
$$\min_{s' }|\theta-\theta_s(a, \tilde\omega)|>e^{-{\bar N}^\kappa},$$
where $\theta_{s'}$ is a family of Lipschitz functions satisfying $\Vert \theta_{s'}\Vert_{\text Lip}\leq C\bar N$, 
$s'<S'$ and $\log\log S'\sim\log\log\bar N$; 
$u=u^{(r)}$, $|\log\delta|^s\leq \bar N\leq r^C$, $s$, $C>1$, $r\geq R$. In fact the above estimates hold for all $u= u^{(r')}$ with $r'>r\geq R$ fulfilling 
assumptions (Hi-iii) and verifying (5.11).
\endproclaim

\demo{Proof}
The first paragraph of the proof of Lemma 5.3 together with the induction 
scheme of Lemma 19.13 in \cite {B3} produce the family of Lipschitz functions $\theta_{s'}$.
The initial family of Lipschitz functions here is just the set of roots (in $\theta$) of the determinants of the various at most 
$(2b+2d)\times (2b+2d)$ matrices.
\hfill $\square$
\enddemo

\proclaim{Lemma 5.5} Assume that Lemmas 5.3 and 5.4 are available, then conditions (Hi-iv) with the bound (5.11) are met for $r=R+1$.  Moreover for all $r\geq 1$, we have the 
bounds:
$$\delta_r<\delta^2 M^{-(\frac{4}{3})^r}, \, \bar\delta_r<\delta^2 M^{-\frac{1}{2}(\frac{4}{3})^r}; \kappa_r<\delta^{2p+4} M^{-(\frac{4}{3})^{r+2}}, \, \bar\kappa_r<\delta^{2p+4} M^{-\frac{1}{2}(\frac{4}{3})^{r+2}}.$$
\endproclaim
\demo{Proof} 
We assume that conditions (Hi-iv) hold at stage $R$. To construct $u^{(R+1)}$, the key is to control $T^{-1}_N(u^{(R)})$ with $N=M^{R+1}$
after a further excision of the $(a,\tilde\omega)$ set. This will produce (Hiv, c) at stage $R+1$. It is as in (19.76-19.86) in \cite{B3}. 
Below we repeat some of the details.

We cover each copy of $[-M^{R+1}, M^{R+1}]^{b+d}$ with $[-M^{R}, M^{R}]^{b+d}$ and cubes $J$ of the form
$$[-L, L]^{b+d}+k,$$ where $L=(\log N)^{C_6}=\Cal O(R+1)^{C_6}$ as in \cite{B3} just above (19.82), $k\in\Bbb Z^{b+d}$ satisfying
$$M^R/2<|k|<M^{R+1}$$ and use the resolvent
identity. 

We first estimate $T^{-1}_{M^R}(u^{(R)})$. Fix $(a, \tilde\omega)\in\bigcup_{I\in\Lambda_R}I$. Condition (Hiv, c) at
stage $R$ gives
$$\align
\Vert T^{-1}_{M^R}(u^{(R-1)})\Vert&< M^{R^C},\tag 5.14\\
|T^{-1}_{M^R}(u^{(R-1)})(\xi,\xi')|&<e^{-|\xi-\xi'|^c}\text { for }|\xi-\xi'|>CR^{C/c}.\tag 5.15\endalign$$
We write
$$T_{M^R}(u^{(R)})=T_{M^R}(u^{(R-1)})+[T_{M^R}(u^{(R)})-T_{M^R}(u^{(R-1)})].$$
Using the resolvent equation, (5.14), condition (Hii) and the bound on $\delta_R$, this gives
$$\Vert T^{-1}_{M^R}(u^{(R)})\Vert
\leq M^{R^C}+ \Cal O(1) M^{-(\frac{4}{3})^R}M^{2R^C}
< 2M^{R^C}.\tag 5.16$$
Using the resolvent series, the above norm bound, condition (iii) at stage $R$ and (5.15) yields 
$$|T^{-1}_{M^R}(u^{(R)})(\xi,\xi')|<e^{-|\xi-\xi'|^c}\text { for } |\xi-\xi'|>CR^{C/c}.\tag 5.17$$

We now study $T(u^{(R)})$ restricted to the $J$ cubes, $T_J(u^{(R)})$. We distinguish two types of 
$J$ cubes in $[-M^{R+1}, M^{R+1}]^{b+d}$: 
\item{(a)} For all $(n, j)\in J$, $|j|\geq L$;
\item{(b)} There exists $(n, j)\in J$, $|j|<L$.  

To control the $J$ cubes in (a), we use Lemma 5.4 here and make excisions as in (19.84) in \cite{B3}
with the measure estimates $$M^{-(R+1)}\ll[\exp\exp(\log (R+1))^{1/3}]^{-1}$$
just below (19.84). The last estimate is of the order of the bound in (Hiv, d). 

The $J$ cubes in (b) are controlled by using Lemma 5.3 here and Lemma 9.9 in\cite{B3}. This is the same as in the paragraph containing (19.85, 19.86) 
combined with the construction starting from (19.76) leading to (19.81) 
in addition to that in (18.28-18.33) in \cite{B3}.  
For the measure estimates we use that the amplitude-frequency map is a diffeomorphism satisfying
$$\Vert \partial\tilde\omega/\partial a\Vert\asymp \Vert( \partial\tilde\omega/\partial a)^{-1}\Vert\asymp 1$$
from (5.8-5.10) and Proposition 4.3.

The conclusion is that 
$$\align
\Vert T_J^{-1}(u^{(R)})\Vert&<e^{L^\sigma},\\
|T^{-1}_{J}(u^{(R)})(\xi,\xi')|&<e^{-|\xi-\xi'|^c}\text { for }|\xi-\xi'|>L/10,\tag 5.18\endalign$$
with $0<\sigma<c<1$, for all $J=[-L, L]^{b+d}+k,$ and $k$ satisfying
$$M^R/2<|k|<M^{R+1}.$$ 

Using (5.16-5.18) and an application of Lemma 5.1 in \cite{BW} adjusted  to the sub-exponential setting, we obtain the resolvent estimate (Hiv, c) at scale $M^{R+1}$. 
We then construct $u^{(R+1)}$ as in \cite{B3, Chap. 18, IV. P140-141} (18.36-18.41), cf. sect. 6 of \cite{BW}. 
We note that as before using the derivative estimates $\partial u^{(R+1)}$, $u^{(R+1)}$ {\it as a function} is defined on the
entire $(a,\tilde \omega)$ space. Therefore ${\tilde\omega}^{(R+1)}$ is obtained by application of the implicit function theorem to the $Q$-equations
and is also defined on the entire $(a,\tilde \omega)$ space.  So (Hi-iv) are available at stage $R+1$.
The induction from $R+1\to$ $R+2$ 
proceeds with $u^{(R+1)}$ replacing $u^{(R)}$ ...
we have therefore proved the Lemma. 
\hfill$\square$
\enddemo
\demo{Proof of the Theorem}
Lemma 5.5 together with (5.9-5.10) prove the Theorem with $a\to a/\delta$, taking into account the rescaling starting in sect. 3, and $\omega_k=j^2_k+\delta^{2p}\tilde\omega_k$ for $k=1, 2, ..., b$. 
\hfill$\square$
\enddemo

\enddocument
\end